\documentclass[11pt,reqno,psamsfonts]{amsart}
\usepackage{amssymb,epsfig}

\allowdisplaybreaks[4]
\numberwithin{equation}{section}
\swapnumbers

\theoremstyle{plain}
\newtheorem{thm}{Theorem}[section]
\newtheorem{lemma}[thm]{Lemma}
\newtheorem{prop}[thm]{Proposition}
\newtheorem{coro}[thm]{Corollary}

\theoremstyle{definition}
\newtheorem{notn}[thm]{Notation}

\newtheorem{remks}[thm]{Remarks}

\newtheorem{partitions}[thm]{Partitions}
\newtheorem{graphs}[thm]{Graphs}
\newtheorem{word-graphs}[thm]{Word-graphs}
\newtheorem{quotient graphs}[thm]{Quotient graphs}
\newtheorem{admissible graphs}[thm]{Admissible graphs and congruences}
\newtheorem{loop-characteristic}[thm]{The loop-characteristic of admissible graphs}
\newtheorem{strongly admissible graphs}[thm]{Strongly admissible graphs and strong congruences}
\newtheorem{permutations}[thm]{Permutations with restricted cycle lengths}
\newtheorem{cycles}[thm]{Cycles of a given length}
\newtheorem{permutations and graphs}[thm]{Permutations compatible with graphs}
\newtheorem{non-commutative probability spaces}[thm]{Non-commutative $*$-probability spaces}
\newtheorem{random matrices}[thm]{Random matrices and convergence in $*$-distribution}
\newtheorem{random permutations}[thm]{Random permutation matrices}
\newtheorem{gaussian matrices}[thm]{Asymptotic freeness from Gaussian matrices}

\newcommand{\A}{\ensuremath{\mathcal A}}
\newcommand{\C}{\ensuremath{\mathbb C}}

\newcommand{\Con}{\ensuremath{\mathrm{Con}}}
\newcommand{\card}{\ensuremath{\mathrm{card}}}
\newcommand{\converges}{\ensuremath{\xrightarrow[k\to\infty]{\text{$*$-distr.}}}}
\newcommand{\convergesas}{\ensuremath{\xrightarrow[k\to\infty]{\text{$*$-distr. a.s.}}}}
\newcommand{\dif}{\ensuremath{\,\mathrm d}}
\newcommand{\E}{\ensuremath{\mathbb E}}
\newcommand{\F}{\ensuremath{\mathbb F}}
\newcommand{\M}{\ensuremath{\mathcal M}}

\newcommand{\N}{\ensuremath{\mathbb N}}
\newcommand{\OFP}{\ensuremath{(\Omega,\mathcal F,P)}}
\newcommand{\on}{\ensuremath{\mathrm O(1/N_k)}}
\newcommand{\onn}{\ensuremath{\mathrm O\left(1/N_k^{1+\delta}\right)}}

\renewcommand{\S}{\ensuremath{\mathcal S}}
\newcommand{\SCon}{\ensuremath{\mathrm{SCon}}}
\newcommand{\toA}{\ensuremath{^{(A)}}}
\newcommand{\tok}{\ensuremath{^{(k)}}}
\newcommand{\toMN}{\ensuremath{^{(M,N)}}}

\newcommand{\tr}{\ensuremath{\mathrm{tr}}}

\begin{document}

\title[Asymptotic freeness of random permutations with restricted cycle lengths]{Asymptotic freeness of random permutation matrices with restricted cycle lengths}
\author{Mihail G. Neagu}
\address{Department of Mathematics \& Statistics, Queen's University, Jeffery Hall, University Avenue, Kingston, Ontario, Canada\\K7L 3N6}
\email{mgneagu@mast.queensu.ca}
\date{March 2007}
\begin{abstract}
Let $A_1,A_2,\ldots,A_s$ be a finite sequence of (not necessarily disjoint, or even distinct) non-empty sets of positive integers such that each $A_r$ either is a finite set or satisfies $\sum\limits_{j\in\N\setminus A_r}\frac{1}{j}<\infty$. It is shown that an independent family $U_1,U_2,\ldots,U_s$ of uniformly distributed random $N\times N$ permutation matrices with cycle lengths restricted to $A_1,A_2,\ldots,A_s$, respectively, converges in $*$-distribution as $N\to\infty$ to a $*$-free family $u_1,u_2,\ldots,u_s$ of non-commutative random variables, where each $u_r$ is a $(\max A_r)$-Haar unitary (if 
$A_r$ is a finite set) or a Haar unitary (if $A_r$ is an infinite set). Under the additional assumption that each of the sets $A_1,A_2,\ldots,A_s$ either consists of a single positive integer or is infinite, it is shown that the convergence in $*$-distribution actually holds almost surely.
\end{abstract}
\maketitle
\section{Introduction.}
Almost from the initial development of free probability theory, it was shown by Voiculescu that, in many important cases, independent families of $N\times N$ random matrices are asymptotically free in the limit $N\to\infty$; as shown in \cite{Voiculescu}, this asymptotic freeness phenomenon occurs for square matrices with independent complex Gaussian entries and for random unitary matrices.

Another instance of asymptotic freeness which arises in a natural situation is the case of an independent family of uniformly distributed $N\times N$ random permutation matrices, which was shown in \cite{Nica1} to provide an asymptotic model as $N\to\infty$ for a $*$-free family of Haar unitaries in a non-commutative probability space. It was shown in \cite{Neagu} that these matrices are also asymptotically free from various types of independent families of complex Gaussian matrices.

Closely related to the distribution of a Haar unitary is that of a $d$-Haar unitary for a positive integer $d$; a $d$-Haar unitary in a $*$-probability space $(\mathcal A,\varphi)$ is a unitary $u\in\mathcal A$ satisfying $u^d=1$ and $\varphi(u^k)=0$ if $k$ is not a multiple of $d$. It is immediate that this distribution is the same as that of an $N\times N$ random permutation matrix which is uniformly distributed over the set
\begin{equation}
\S_N^{(d)}:=\{\sigma\in\S_N\,|\,\text{all cycles of $\sigma$ have length $d$}\}
\end{equation}
(with $N$ being a multiple of $d$, so that the set $\S_N^{(d)}$ is non-empty). A natural question is then whether asymptotic freeness still occurs for an independent family $U_1,U_2,\ldots,U_s$, where each $U_r$ is an $N\times N$ random permutation matrix which is uniformly distributed over either $\S_N$ or $\S_N^{(d_r)}$ for some positive integer $d_r$ (with $N$ chosen such that each of the sets $\S_N^{(d_r)}$ is non-empty). This paper gives an affirmative answer to a more general form of this question, and provides a unified framework for the two cases in which $U_r$ is uniformly distributed over $\S_N$ or  $\S_N^{(d_r)}$.

More concretely, for every non-empty set $A$ of positive integers, let
\begin{equation}
\S_N\toA:=\{\sigma\in\S_N\,|\,\text{the length of every cycle of $\sigma$ belongs to $A$}\};
\end{equation}
the cases in which $A$ consists either of all positive integers or of a single positive integer were mentioned above. Let $u_1,u_2,\ldots,u_s$ be a $*$-free family of random variables in a $*$-probability space, with each $u_r$ being either a Haar unitary or a $d_r$-Haar unitary for some positive integer $d_r$. Then the question asked in the previous paragraph is whether an asymptotic model for the family $u_1,u_2,\ldots,u_s$ is given by an independent family $U_1,U_2,\ldots,U_s$, where each $U_r$ is an $N\times N$ random permutation matrix which is uniformly distributed over $\S_N^{(A_r)}$, with
\begin{equation}
A_r:=\begin{cases}\{d_r\}&\text{if $u_r$ is a $d_r$-Haar unitary}
\\\N&\text{if $u_r$ is a Haar unitary}
\end{cases}
\end{equation}
(and $N$ chosen such that each of the sets $\S_N^{(A_r)}$ is non-empty). The main result of this paper (Theorem \ref{thm} below) is that the above statement is true and that, in fact, it can be considerably strengthened by allowing 
\begin{equation}\label{ar}
A_r=\begin{cases}\text{any finite set of positive integers with $d_r=\max A_r$}&\text{if $u_r$ is a $d_r$-Haar unitary}
\\\text{any infinite set of positive integers with $\sum\limits_{j\in\N\setminus A_r}\frac{1}{j}<\infty$}&\text{if $u_r$ is a Haar unitary.}
\end{cases}
\end{equation}
Under the additional assumption that each of the sets $A_r$ either consists of a single positive integer or is infinite, it is shown that the convergence in $*$-distribution of the family $U_1,U_2,\ldots,U_s$ actually holds almost surely.

Moreover, the same techniques used in the proof of Theorem 3.1 of \cite{Neagu} can be applied to show that the family $U_1,U_2,\ldots,U_s$ is asymptotically $*$-free from various types of independent families of complex Gaussian matrices, with the convergence in $*$-distribution holding almost surely under the same additional assumption as in the previous paragraph (see \ref{gaussian matrices} below).

The proof of the main theorem of the paper makes use of two principal ideas, which are explored in Sections 2 and 3. The first involves the analysis of a certain class of directed, edge-colored graphs and their congruences, which were introduced in \cite{Nica2}. Section 2 of the paper is devoted to the study of the relevant facts about these graphs, and an important ingredient is a quantity which bears an uncanny resemblance to the classical Euler characteristic. Because of this resemblance, I have called this quantity the \emph{loop-characteristic} of a graph (see \ref{loop-characteristic} below); as far as I have been able to determine, it does not appear in the graph theory literature.

Section 3 of the paper is devoted to the evaluation of certain probabilities related to random permutations which are uniformly distributed over $\S_N\toA$. More specifically, for a fixed non-empty set $A$ of positive integers and for a fixed $k\in A$, the relevant concept for the present purposes is the asymptotic behavior as $N\to\infty$ (with $\S_N\toA\neq\emptyset$) of the expected number of cycles of length $k$ of such a random permutation. If $A$ is a finite set, this asymptotic behavior can be computed with a moderate amount of effort by using Hayman's method (a saddle point method) for determining the asymptotic behavior of power series coefficients of certain complex analytic functions; here Hayman's method is applicable to the exponential generating function for the number of permutations in $\S_N\toA$. If $A$ is an infinite set, it turns out that this asymptotic behavior has connections to recent research (see \cite{Hildebrand} and \cite{Manstavicius}), which shows that it can be computed under the assumption that $A$ satisfies the condition appearing in the second part of \eqref{ar}.

Finally, Section 4 of the paper sets the free probabilistic framework which is necessary for the formulation of the main result, which is then stated and proved.

I am deeply grateful to my former thesis advisor, Professor Alexandru Nica, for his support during the writing of this work.
\section{Edge-Colored Graphs and Their Congruences.}
This section discusses a class of graphs which is needed in the sequel. These graphs were introduced in \cite{Nica2}, which contains additional information about their congruences and homomorphisms. With slight modifications, the definitions and notation appearing in \ref{graphs}, \ref{words}, \ref{quotient graphs}, and  \ref{admissible graphs} below are taken from that work.
\begin{notn}\label{notn s}
(a) Throughout the paper, $s$ is a fixed positive integer and $d_1,d_2,\ldots,d_s$ are fixed elements of the set $\{1,2,\ldots\}\cup\{\infty\}$. If $n$ is a positive integer, then $[n]$ denotes the set
\begin{equation}
[n]:=\{1,2,\ldots,n\}.
\end{equation}

(b) The free monoid generated by some fixed symbols $g_1,g_1^*,g_2,g_2^*,\ldots,g_s,g_s^*$ is denoted by \F. The elements of \F\ are ``words'' in $g_1,g_1^*,g_2,g_2^*,\ldots,g_s,g_s^*$. The identity of \F\ is the ``empty word'' and is denoted by $e$. 

We shall denote by $\approx$ the congruence of $\F$ generated by the relations
\begin{equation}\label{approx}
\big\{g_rg_r^*\approx g_r^*g_r\approx e\,\big|\,r\in[s]\big\}\cup\big\{g_r^{d_r}\approx(g_r^*)^{d_r}\approx e\,\big|\,r\in[s];d_r<\infty\big\}.
\end{equation}
In other words, $\approx$ is the smallest equivalence relation on $\F$ containing the above relations and having the property  that  $w_1\approx w_2\implies uw_1v\approx uw_2v$ for every $u,v,w_1,w_2\in\F$. It is clear that for every $w\in\F\setminus\{e\}$, we have that $w\approx e$ iff some cyclic permutation of the symbols of $w$ is of the form $vx$, where
$v\approx e$ and
\begin{equation}\label{x}
x\in\big\{g_rg_r^*,g_r^*g_r\,\big|\,r\in[s]\big\}\cup\big\{g_r^{d_r},(g_r^*)^{d_r}\,\big|\,r\in[s];d_r<\infty\big\}.
\end{equation}
\end{notn}
\begin{partitions}\label{partitions}
Let $V$ be a non-empty set.

(a) The set of all partitions of $V$ is denoted by $\Pi(V)$. If $\pi\in\Pi(V)$ and $a\in V$, then $\pi(a)$ denotes the block of $\pi$ to which $a$ belongs, $|\pi|$ denotes the number of blocks of $\pi$, and $\overset{\pi}{\sim}$ denotes the equivalence relation on $V$ associated to $\pi$:
\begin{equation}
a\overset{\pi}{\sim}b\overset{\text{def}}{\iff}\pi(a)=\pi(b)\qquad\forall a,b\in V.
\end{equation}
If $\pi\in\Pi(V)$ and $W$ is a non-empty subset of $V$, then the \emph{restriction} of $\pi$ to $W$ is the partition $\rho\in\Pi(W)$ defined by
\begin{equation}
a\overset{\rho}{\sim}b\overset{\text{def}}{\iff}a\overset{\pi}{\sim}b\qquad\forall a,b\in W.
\end{equation}
\pagebreak

(b) For every $\pi,\rho\in\Pi(V)$, we write $\pi\leq\rho$ iff every block of $\pi$ is completely contained in a block of $\rho$. This gives a partial order which makes $\Pi(V)$ into a complete lattice, the meet operation of which is denoted by $\wedge$. If $\pi,\rho\in\Pi(V)$ and $\pi\leq\rho$, then $\rho/\pi$ denotes the partition of $\pi$ defined by
\begin{equation}\label{quotient partitions}
\pi(a)\overset{\rho/\pi}{\sim}\pi(b)\overset{\text{def}}{\iff}a\overset{\rho}{\sim}b\quad\forall a,b\in V.
\end{equation}
\end{partitions}
\begin{graphs}\label{graphs}
(a) Throughout the paper, \emph{graph} means a finite, directed, edge-colored graph with colors from the set $[s]$. Formally, a graph $\Gamma$ consists of a non-empty finite \emph{vertex set} $V$ together with (not necessarily disjoint) subsets $E_1,E_2,\ldots,E_s$ of $V\times V$; for every $r\in[s]$ and for every $a,b\in V$, we shall say that $a\rightarrow b$ is an \emph{$r$-edge} (or an \emph{edge of color $r$}) of $\Gamma$ iff $(a,b)\in E_r$. A graph is said to be \emph{trivial} iff it has no edges (of any color). 

A graph is said to be \emph{monocolored} iff all its edges have the same color. For every $r\in[s]$, the (monocolored) graph obtained from a graph $\Gamma$ by removing all edges which are not of color $r$ is denoted by $\Gamma(r)$.

(b) A \emph{subgraph} of a graph $\Gamma$ is a graph $\Gamma'$ which has vertex set contained in that of $\Gamma$ and which has the property that for every $r\in[s]$, each $r$-edge of $\Gamma'$ is an $r$-edge of $\Gamma$.

A graph $\Gamma$ is said to be \emph{connected} iff for every two distinct vertices $a$ and $b$ of $\Gamma$, there exist a positive integer $n$ and vertices $a_0,a_1,\ldots,a_n$ of $\Gamma$ with $a_0=a$ and $a_n=b$ such that for every $k\in[n]$, there exists an edge (of any color and either direction) between $a_{k-1}$ and $a_k$. A \emph{connected component} of a graph $\Gamma$ is a maximal connected subgraph of $\Gamma$.

(c) Two graphs $\Gamma_1$ and $\Gamma_2$ with vertex sets $V_1$ and $V_2$, respectively, are said to \emph{isomorphic}, written $\Gamma_1\cong\Gamma_2$, iff there exists a bijection $\Phi:V_1\to V_2$ (called an \emph{isomorphism}) such that for every $r\in[s]$ and for every $a,b\in V_1$, we have that
$$\text{$\Phi(a)\rightarrow\Phi(b)$ is an $r$-edge of $\Gamma_2$}\iff\text{$a\rightarrow b$ is an $r$-edge of $\Gamma_1$}.$$
\end{graphs}
\begin{word-graphs}\label{words}
Let \F\ be as in \ref{notn s}(b). For every $w:=g_{r_1}^{\varepsilon_1}g_{r_2}^{\varepsilon_2}\cdots g_{r_n}^{\varepsilon_n}\in\F\setminus\{e\}$, where $n$ is a positive integer, $r_1,r_2,\ldots,r_n\in[s]$, and $\varepsilon_1,\varepsilon_2,\ldots,\varepsilon_n\in\{1,*\}$, the \emph{word-graph} $\Gamma_w$ associated to $w$ is the graph with vertex set $[n]$ and with $n$ edges, described as follows: for every $k\in[n]$, the $k$-th edge of $\Gamma_w$ has color $r_k$ and is
\begin{equation}
\begin{cases}k\;\overset{r_k}{\longleftarrow}\;k+1&\text{if $\varepsilon_k=1$}
\\k\;\overset{r_k}{\longrightarrow}\;k+1&\text{if $\varepsilon_k=*$}
\end{cases}
\end{equation}
(with the convention that $n+1:=1$). Figure \ref{fig1} below shows an example of a word-graph.
\end{word-graphs}
\begin{figure}[ht]
\centerline{\input{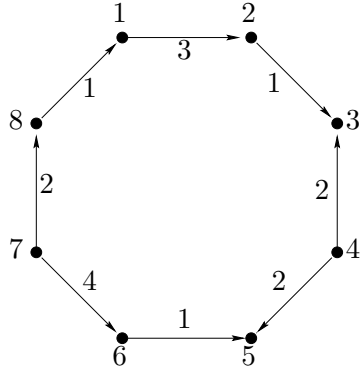}}
\caption{\label{fig1} The word-graph associated to the word $g_3^*g_1^*g_2g_2^*g_1g_4g_2^*g_1^*$.}
\end{figure}
\begin{quotient graphs}\label{quotient graphs}
For every graph $\Gamma$ and for every partition $\pi$ of the vertex set of $\Gamma$, the \emph{quotient} of $\Gamma$ by $\pi$ is the graph $\Gamma/\pi$ which has vertex set $\pi$ and which is defined by requiring that for every $r\in[s]$ and for every pair of vertices $a$ and $b$ of $\Gamma$, we have that $\pi(a)\rightarrow\pi(b)$ is an $r$-edge of $\Gamma/\pi$ iff there exists an $r$-edge $a'\rightarrow b'$ of $\Gamma$ with $a'\overset{\pi}{\sim}a$ and $b'\overset{\pi}{\sim}b$. Figure \ref{fig2} below shows an example of a quotient graph.

It is immediate that any quotient of a connected graph is connected; in particular, with the notation from \ref{words}, we have that $\Gamma_w/\pi$ is connected for every $w\in\F\setminus\{e\}$ and for every partition $\pi$ of the vertex set of $\Gamma_w$.
\end{quotient graphs}
\begin{figure}[ht]
\centerline{\input{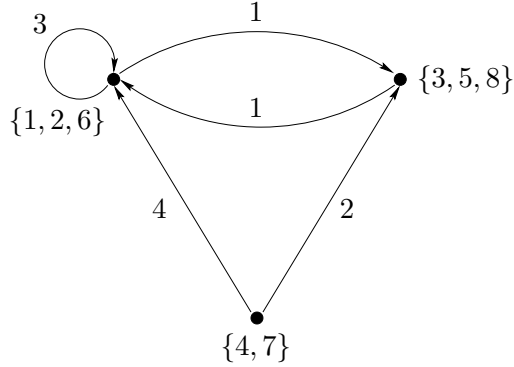}}
\caption{\label{fig2} The quotient of the graph from Figure \ref{fig1} by the partition
$\big\{\{1,2,6\},\{3,5,8\},\{4,7\}\big\}$.}
\end{figure}
\begin{lemma}[Second Isomorphism Theorem]\label{second isomorphism theorem}
If $\Gamma$ is a graph with vertex set $V$ and $\pi,\rho\in\Pi(V)$ are such that $\pi\leq\rho$, then $\Gamma/\rho\cong(\Gamma/\pi)/(\rho/\pi)$, where $\rho/\pi$ is as in \ref{partitions}(b).
\begin{proof}
Let $\Phi:\rho\to\rho/\pi$ be defined by 
$$\Phi\big(\rho(a)\big)=(\rho/\pi)\big(\pi(a)\big)\qquad\forall a\in V.$$ 
From the definition of $\rho/\pi$, it is immediate  that $\Phi$ is a well-defined bijection. For every $r\in[s]$ and for every $a,b\in V$, we have that
\begin{align*}
(\rho/\pi)&\big(\pi(a)\big)\rightarrow(\rho/\pi)\big(\pi(b)\big)\text{is an $r$-edge of $(\Gamma/\pi)/(\rho/\pi)$}
\\&\iff\text{there exists an $r$-edge $A\rightarrow B$ of $\Gamma/\pi$ such that $A\overset{\rho/\pi}{\sim}\pi(a)$ and $B\overset{\rho/\pi}{\sim}\pi(b)$}
\\&\iff\text{there exists an $r$-edge $a'\rightarrow b'$ of $\Gamma$ such that $\pi(a')\overset{\rho/\pi}{\sim}\pi(a)$ and $\pi(b')\overset{\rho/\pi}{\sim}\pi(b)$}
\\&\iff\text{there exists an $r$-edge $a'\rightarrow b'$ of $\Gamma$ such that $a'\overset{\rho}{\sim}a$ and $b'\overset{\rho}{\sim}b$}
\\&\iff\text{$\rho(a)\rightarrow\rho(b)$ is an $r$-edge of $\Gamma/\rho$.}
\end{align*}
Hence $\Phi$ is an isomorphism between $\Gamma/\rho$ and $(\Gamma/\pi)/(\rho/\pi)$.
\end{proof}
\end{lemma}
\begin{admissible graphs}\label{admissible graphs}
(a) A graph $\Gamma$ is said to be \emph{admissible} iff no two distinct edges of $\Gamma$ having the same color begin at the same vertex or end at the same vertex.

(b) A \emph{congruence} of a graph $\Gamma$ is a partition $\pi$ of the vertex set of $\Gamma$ such that for every pair $a_1\rightarrow b_1$ and $a_2\rightarrow b_2$ of edges of $\Gamma$ having the same color, we have that $a_1\overset{\pi}{\sim}a_2\iff b_1\overset{\pi}{\sim}b_2$. The set of all congruences of $\Gamma$ is denoted by $\Con(\Gamma)$. It is easily seen that
\begin{equation}
\pi\wedge\rho\in\Con(\Gamma)\qquad\forall\pi,\rho\in\Con(\Gamma)
\end{equation}
and that for every partition $\pi$ of the vertex set of $\Gamma$, we have that
\begin{equation}\label{admissible}
\Gamma/\pi\text{ is admissible }\iff\pi\in\Con(\Gamma).
\end{equation}
(c) Let $\Gamma$ be an admissible graph. For every $r\in[s]$, an \emph{$r$-path} (or a \emph{path of color $r$}) of
$\Gamma$ is a non-trivial connected component of $\Gamma(r)$ (notation and terminology as in \ref{graphs}). The \emph{length} $l(P)$ of a path $P$ (of any color) of $\Gamma$ is the number of edges of $P$. The graph obtained from $\Gamma$ by
removing all edges of a path $P$ is denoted by $\Gamma\setminus P$.
\end{admissible graphs}
We shall make use of the following observation.
\begin{lemma}\label{lemma paths}
If $\Gamma$ is a non-trivial connected admissible graph, then there exists a path $P$ of $\Gamma$ such that $\Gamma\setminus P$ has at most one non-trivial connected component.
\begin{proof}
The proof is by induction on the number of paths of $\Gamma$. If $\Gamma$ consists of a single path $P$, then it is clear that $\Gamma\setminus P$ has no non-trivial connected components.

Now suppose that the assertion of the lemma holds for every non-trivial connected admissible graph with number of paths less than or equal to some positive integer $n$, and let $\Gamma$ be a non-trivial connected admissible graph with $n+1$ paths. Let $\Gamma'$ be a subgraph of $\Gamma$ with minimal number of paths relative to the following property: there exists a path $Q$ of $\Gamma$ such that $\Gamma'$ is a non-trivial connected component of $\Gamma\setminus Q$.

If $\Gamma'$ consists of a single path $P$, then it is clear that $\Gamma\setminus P$ has exactly one non-trivial connected component. Finally, suppose that $\Gamma'$ has at least two paths. By the induction hypothesis, there exists a path $P$ of $\Gamma'$ such that $\Gamma'\setminus P$ has exactly one non-trivial connected component $\Gamma''$. If $\Gamma''$ has no vertices in common with $Q$, then $\Gamma''$ is a non-trivial connected component of $\Gamma\setminus P$, which contradicts the minimality of the number of paths of $\Gamma'$. Thus $\Gamma''$ has at least one vertex in common with $Q$, from which it follows that $\Gamma\setminus P$ has exactly one non-trivial connected component.
\end{proof}
\end{lemma}
\begin{loop-characteristic}\label{loop-characteristic}
Let $\Gamma$ be an admissible graph. It is clear that every path of $\Gamma$ is a circuit which is not self-intersecting. For every $r\in[s]$, an \emph{$r$-loop} (or a \emph{loop of color $r$}) of $\Gamma$ is an $r$-path of $\Gamma$ which is a closed circuit, and an \emph{$r$-string} (or a \emph{string of color $r$}) of $\Gamma$ is an $r$-path of $\Gamma$ which is not a loop. The \emph{loop-characteristic} of $\Gamma$ is the integer $\chi(\Gamma)$ defined by
\begin{equation}
\chi(\Gamma)=V(\Gamma)-E(\Gamma)+L(\Gamma),
\end{equation}
where $V(\Gamma)$, $E(\Gamma)$, and $L(\Gamma)$ denote the numbers of vertices, edges, and loops (of any color) of $\Gamma$, respectively. For instance, the (admissible) graph shown in Figure \ref{fig2} above has a 1-loop of length 2, a 2-string of length 1, a 3-loop of length 1, a 4-string of length 1, and has loop-characteristic $3-5+2=0$.
\end{loop-characteristic}
\begin{prop}\label{connected chi}
If $\Gamma$ is a connected admissible graph, then $\chi(\Gamma)\leq1$, with equality holding iff each path $P$ of $\Gamma$ has the following property: every connected component $C$ of $\Gamma\setminus P$ has exactly one vertex in common with $P$ and satisfies $\chi(C)=1$. Moreover, if some path of $\Gamma$ has this property, then $\chi(\Gamma)=1$.
\begin{proof}
The proof is by induction on the number of paths of $\Gamma$. If $\Gamma$ has no paths (and therefore no edges and exactly one vertex), then it is immediate that $\chi(\Gamma)=1$.

Now suppose that the assertion of the proposition holds for every connected admissible graph with number of paths less than some positive integer $n$, let $\Gamma$ be a connected admissible graph with $n$ paths, and let $P$ be any path of $\Gamma$. Let $C_1,C_2,\ldots,C_k$ be the connected components of $\Gamma\setminus P$ and let $n_i$ be the number of vertices which $P$ has in common with each $C_i$. Then
\begin{align*}
\chi(\Gamma)&=\chi(P)+\chi(C_1)+\chi(C_2)+\cdots+\chi(C_k)-(n_1+n_2+\cdots+n_k)\tag*{\text{[by the definition of $\chi$]}}
\\&\leq1+k-(n_1+n_2+\cdots+n_k)\tag*{\text{[since $\chi(P)=1$ and by the induction hypothesis]}}
\\&\leq1,\tag*{\text{[since $\Gamma$ is connected, so $n_1,n_2,\ldots,n_k\geq1$]}}
\end{align*}
with equalities holding iff $\chi(C_1)=\chi(C_2)=\cdots=\chi(C_k)=1$ and $n_1=n_2=\cdots=n_k=1$.
\end{proof}
\end{prop}
\begin{strongly admissible graphs}\label{strongly admissible graphs}
(a) A graph $\Gamma$ is said to be \emph{strongly admissible} iff it is admissible and for every $r\in[s]$, each $r$-loop of $\Gamma$ has length $d_r$ and each $r$-string of $\Gamma$ has length strictly less than $d_r$ 
[notation as in \ref{notn s}(a)]. In particular, a strongly admissible graph cannot contain any $r$-loops 
for any $r\in[s]$ with $d_r=\infty$.

(b) A \emph{strong congruence} of a graph $\Gamma$ is a congruence $\pi$ of $\Gamma$ such that $\Gamma/\pi$ is strongly admissible. The set of all strong congruences of $\Gamma$ is denoted by $\SCon(\Gamma)$. It is easily seen that
\begin{equation}\pi\wedge\rho\in\SCon(\Gamma)\qquad
\forall\pi,\rho\in\SCon(\Gamma).
\end{equation}
\end{strongly admissible graphs}
\begin{remks}\label{subgraphs}
Let $\Gamma'$ be a subgraph of a graph $\Gamma$.
\vspace{2mm}

(a) If $\Gamma$ is admissible (respectively, strongly admissible), then so is $\Gamma'$.
\vspace{2mm}

(b) If both $\Gamma$ and $\Gamma'$ are connected and admissible and $\chi(\Gamma)=1$, then $\chi(\Gamma')=1$.
\vspace{2mm}

Indeed, (a) is immediate from the definition of admissible and strongly admissible graphs, whereas (b) follows from Proposition \ref{connected chi}, since each path $P'$ of $\Gamma'$ is a subgraph of some path $P$ of $\Gamma$ and every connected component of $\Gamma'\setminus P'$ is a subgraph of some connected component of $\Gamma\setminus P$.
\end{remks}
\begin{lemma}\label{lemma scon chi}
If $\Gamma$ is a connected strongly admissible graph and $\pi\in\SCon(\Gamma)$ is such that $\chi(\Gamma/\pi)=1$, then all blocks of $\pi$ are singletons.
\begin{proof}
The proof is by induction on the number of paths of $\Gamma$. If $\Gamma$ has no paths (and therefore no edges and exactly one vertex), then the assertion holds vacuously, since the only partition of the vertex set of $\Gamma$ consists of a
singleton block.

Now suppose that the assertion of the lemma holds for every connected strongly admissible graph with number of paths less than some positive integer $n$, let $\Gamma$ be a connected strongly admissible graph with $n$ paths, and let
$\pi\in\SCon(\Gamma)$ be such that $\chi(\Gamma/\pi)=1$. Let $P$ be any path of $\Gamma$ (of some color $r$). Since $\Gamma$ is strongly admissible, $P$ is either a loop of length $d_r$ or a string of length strictly less than $d_r$, which implies that no two distinct vertices of $P$ belong to the same block of $\pi$ (since in that case $\Gamma/\pi$ would have an $r$-loop of length strictly less than $d_r$, contradicting the assumption that $\Gamma/\pi$ is strongly admissible). Let $a$ and $b$ be two distinct vertices of $\Gamma$; it will be shown that $a\not\overset{\pi}{\sim}b$.

If $a$ and $b$ belong to the same connected component $C$ of $\Gamma\setminus P$, let $\rho$ be the restriction of $\pi$ to the vertices of $C$, as in \ref{partitions}(a). Since $C$ is a connected subgraph of $\Gamma$ and 
$C/\rho$ is isomorphic to a connected subgraph of $\Gamma/\pi$, it follows from Remarks \ref{subgraphs} and the induction hypothesis that all blocks of $\rho$ are singletons; in particular, $a\not\overset{\pi}{\sim}b$.

Finally, suppose that $a$ and $b$ belong to distinct connected components $C_1$ and $C_2$ of $\Gamma\setminus P$. Since $\Gamma$ is connected, there exist vertices $c$ and $d$ of $P$ such that $c$ belongs to $C_1$ and $d$ belongs to $C_2$. Let $\overline P$ be the $r$-path of $\Gamma/\pi$ containing the vertices $\pi(c)$ and $\pi(d)$, and let $\overline{C_1}$ and $\overline{C_2}$ be the connected components of $(\Gamma/\pi)\setminus\overline P$ containing the vertices $\pi(a)$ and $\pi(b)$, respectively. Since $\chi(\Gamma/\pi)=1$, it follows from Proposition \ref{connected chi} that $\overline P$ has exactly one vertex in common with $\overline{C_1}$ [namely $\pi(c)$] and with $\overline{C_2}$ [namely $\pi(d)$]. Since $\pi(c)\neq\pi(d)$ (as $c$ and $d$ are distinct vertices of $P$), it follows that $\overline{C_1}\neq\overline{C_2}$. In particular, $\pi(a)\neq\pi(b)$, so $a\not\overset{\pi}{\sim}b$.
\end{proof}
\end{lemma}
\begin{prop}\label{prop scon chi}
For every connected graph $\Gamma$, there exists at most one $\pi\in\SCon(\Gamma)$ such that $\chi(\Gamma/\pi)=1$.
\begin{proof}
Let $\pi,\rho\in\SCon(\Gamma)$ be such that $\chi(\Gamma/\pi)=\chi(\Gamma/\rho)=1$. Then $\pi\wedge\rho\in\SCon(\Gamma)$ and
$\Gamma/\pi\cong(\Gamma/\pi\wedge\rho)/(\pi/\pi\wedge\rho)$ by Lemma \ref{second isomorphism theorem}, so all blocks of $\pi/\pi\wedge\rho$ are singletons by Lemma \ref{lemma scon chi}, which implies that $\pi\leq\rho$. A similar argument shows that $\rho\leq\pi$, so $\pi=\rho$.
\end{proof}
\end{prop}
\begin{lemma}\label{lemma scon chi word-graphs}
Let \F\ be as in \ref{notn s}(b), let $v\in\F\setminus\{e\}$, and let $x$ be as in \eqref{x}. Then:
\vspace{2mm}

(a) $\mathrm{card}\big\{\pi\in\SCon(\Gamma_x)\,\big|\,\chi(\Gamma_x/\pi)=1\big\}=1$.
\vspace{2mm}

(b) $\mathrm{card}\big\{\pi\in\SCon(\Gamma_{vx})\,\big|\,\chi(\Gamma_{vx}/\pi)=1\big\}=\mathrm{card}\big\{\rho\in\SCon(\Gamma_v)\,\big|\,\chi(\Gamma_v/\rho)=1\big\}$.
\begin{proof}
Let $r\in[s]$ be such that $x\in\left\{g_rg_r^*,g_r^*g_r,g_r^{d_r},(g_r^*)^{d_r}\right\}$.
\vspace{2mm}

(a) By Proposition \ref{prop scon chi}, it suffices to show that there exists $\pi\in\SCon(\Gamma_x)$ such that $\chi(\Gamma_x/\pi)=1$. If $x=g_rg_r^*$ or $g_r^*g_r$ and $d_r=1$, let $\pi=\big\{\{1,2\}\big\}$. If $x=g_rg_r^*$ or $g_r^*g_r$ and $d_r>1$, let $\pi=\big\{\{1\},\{2\}\big\}$. Finally, if $d_r<\infty$ and $x=g_r^{d_r}$ or $(g_r^*)^{d_r}$, let $\pi=\big\{\{1\},\{2\},\ldots,\{d_r\}\big\}$. In any of the three cases, it is easily seen that $\pi\in\SCon(\Gamma_x)$ and that $\chi(\Gamma_x/\pi)=1$.
\vspace{2mm}

(b) By Proposition \ref{prop scon chi}, it suffices to show that there exists $\pi\in\SCon(\Gamma_{vx})$ such that $\chi(\Gamma_{vx}/\pi)=1$ iff there exists $\rho\in\SCon(\Gamma_v)$ such that $\chi(\Gamma_v/\rho)=1$. Let $n$ be the length of $v$.
\vspace{2mm}

($\Longrightarrow$) Suppose that there exists $\pi\in\SCon(\Gamma_{vx})$ such that $\chi(\Gamma_{vx}/\pi)=1$, and let $\rho$ be the restriction of $\pi$ to $[n]$, as in \ref{partitions}(a).

First note that $1\overset{\pi}{\sim}n+1$; indeed, if $x=g_rg_r^*$ or $g_r^*g_r$, this is because $\pi$ is a congruence, whereas, if $d_r<\infty$ and $x=g_r^{d_r}$ or $(g_r^*)^{d_r}$, it is because the $r$-path of $\Gamma_{vx}/\pi$ which contains $\pi(1)$ must be a loop of length $d_r$. This implies that $\Gamma_v/\rho$ is isomorphic to a connected subgraph of $\Gamma_{vx}/\pi$, so $\rho\in\SCon(\Gamma_v)$ and $\chi(\Gamma_v/\rho)=1$ by Remarks \ref{subgraphs}.
\pagebreak

($\Longleftarrow$) Suppose that there exists $\rho\in\SCon(\Gamma_v)$ such that $\chi(\Gamma_v/\rho)=1$. The two cases in which $x=g_rg_r^*$ or $g_r^*g_r$, and in which $d_r<\infty$ and $x=g_r^{d_r}$ or $(g_r^*)^{d_r}$ are treated separately. In
each of the two cases, there are several subcases to consider, depending on the structure of the graph $\Gamma_v/\rho$.
\vspace{2mm}

\underline{Case 1:} $x=g_rg_r^*$ (the case when $x=g_r^*g_r$ is analogous). 
\vspace{2mm}

If $\Gamma_v/\rho$ has an $r$-edge $A\to\rho(1)$ for some block $A$ of $\rho$, let $\pi$ be the partition of $[n+2]$ which is obtained from $\rho$ by adding $n+1$ to the block $\rho(1)$ and $n+2$ to the block $A$. Then $\Gamma_{vx}/\pi$ is isomorphic to $\Gamma_v/\rho$, so $\pi\in\SCon(\Gamma_{vx})$ and $\chi(\Gamma_{vx}/\pi)=1$.

If $d_r=1$ and $\Gamma_v/\rho$ has no $r$-path containing $\rho(1)$, let $\pi$ be the partition of $[n+2]$ which is obtained from $\rho$ by adding $n+1$ and $n+2$ to the block $\rho(1)$. Then $\Gamma_{vx}/\pi$ is isomorphic to $\Gamma_v/\rho$ with an extra $r$-loop of length 1 added at the vertex $\rho(1)$, so $\pi\in\SCon(\Gamma_{vx})$ and $\chi(\Gamma_{vx}/\pi)=1$.

If $d_r>1$ and $\Gamma_v/\rho$ has no $r$-path containing $\rho(1)$, let $\pi$ be the partition of $[n+2]$ which is obtained from $\rho$ by adding $n+1$ to the block $\rho(1)$ and creating a singleton block $\{n+2\}$. Then $\Gamma_{vx}/\pi$
is isomorphic to $\Gamma_v/\rho$ with an extra $r$-string (of length 1) consisting of the edge $\{n+2\}\to\rho(1)$, so $\pi\in\SCon(\Gamma_{vx})$ and $\chi(\Gamma_{vx}/\pi)=1$.

If $d_r>1$ and $\Gamma_v/\rho$ has an $r$-string $P$ of the form $\rho(1)\to A_1\to A_2\to\cdots A_{d_r-1}$ (of length $d_r-1$), let $\pi$ be the partition of $[n+2]$ which is obtained from $\rho$ by adding $n+1$ to the block 
$\rho(1)$ and $n+2$ to the block $A_{d_r-1}$. Then $\Gamma_{vx}/\pi$ is isomorphic to $\Gamma_v/\rho$ with an extra edge $A_{d_r-1}\to\rho(1)$ (which extends $P$ to a loop of length $d_r$), so $\pi\in\SCon(\Gamma_{vx})$ and $\chi(\Gamma_{vx}/\pi)=1$.

If $d_r>2$ and $\Gamma_v/\rho$ has an $r$-string $P$ of the form $\rho(1)\to A_1\to A_2\to\cdots A_k$ (of length $k$) for some $1\leq k<d_r-1$, let $\pi$ be the partition of $[n+2]$ which is obtained from $\rho$ by adding $n+1$ to the block $\rho(1)$ and creating a singleton block $\{n+2\}$. Then $\Gamma_{vx}/\pi$ is isomorphic to $\Gamma_v/\rho$ with an extra edge $\{n+2\}\to\rho(1)$ (which extends $P$ to a string of length $k+1<d_r$), so $\pi\in\SCon(\Gamma_{vx})$ and $\chi(\Gamma_{vx}/\pi)=1$.
\vspace{2mm}

\underline{Case 2:} $d_r<\infty$ and $x=g_r^{d_r}$ (the case when $x=(g_r^*)^{d_r}$ is analogous).
\vspace{2mm}

If $\Gamma_v/\rho$ has no $r$-path containing $\rho(1)$, let $\pi$ be the partition of $[n+d_r]$ which is obtained from $\rho$ by adding $n+1$ to the block $\rho(1)$ and creating singleton blocks $\{n+2\},\{n+3\},\ldots\{n+d_r\}$. Then  $\Gamma_{vx}/\pi$ is isomorphic to $\Gamma_v/\rho$ with an extra $r$-loop (of length $d_r$) which has no vertices in common with $\Gamma_v/\rho$ other than $\rho(1)$, so $\pi\in\SCon(\Gamma_{vx})$ and $\chi(\Gamma_{vx}/\pi)=1$.

If $\Gamma_v/\rho$ has an $r$-path $P$ of the form $A_j\to\cdots\to A_2\to A_1\to\rho(1)\to B_1\to B_2\to\cdots\to B_k$ (where at least one of $j$ and $k$ is non-zero, $j+k\leq d_r$, and possibly $A_j=B_k$), let $\pi$ be partition of $[n+d_r]$ which is obtained from $\rho$ by adding $n+1$ to the block $\rho(1)$, $n+1+i$ to the block $A_i$ for every $i\in[j]$, $n+d_r+1-i$ to the block $B_i$ for every $i\in[k]$, and (if $j+k<d_r-1$) creating singleton blocks $\{n+j+2\},\{n+j+3\},\ldots,\{n+d_r-k\}$. Then $\Gamma_{vx}/\pi$ is isomorphic to $\Gamma_v/\rho$ with the path $P$ extended to an $r$-loop (of length $d_r$) which has no vertices in common with $\Gamma_v/\rho$ other than the vertices of $P$, so $\pi\in\SCon(\Gamma_{vx})$ and $\chi(\Gamma_{vx}/\pi)=1$.
\end{proof}
\end{lemma}
\pagebreak

\begin{lemma}\label{lemma cyclic}
Let \F\ be as in in \ref{notn s}(b), let $w\in\F\setminus\{e\}$, and let $\pi\in\SCon(\Gamma_w)$. If $\Gamma_w/\pi$ has a path $P$ such that at most one vertex of $P$ belongs to other paths of $\Gamma_w/\pi$, then some cyclic permutation of the symbols of $w$ is of the form $vx$ with $x$ as in \eqref{x}.
\begin{proof}
It may be assumed that $\Gamma_w$ is admissible (otherwise it has a mono-colored subgraph of the form $\bullet\rightarrow\bullet\leftarrow\bullet$ or $\bullet\leftarrow\bullet\rightarrow\bullet$, which immediately implies that some cyclic permutation of the symbols of $w$ is of the form $vu_r^*u_r$ or $vu_ru_r^*$ for some $r\in[s]$). Let $r$ be the color of the path $P$ and let $B$ be a vertex of $P$ (that is, block of $\pi$) such that no other vertex of $P$ belongs to other paths of $\Gamma_w/\pi$. Let $Q$ be an $r$-path of $\Gamma_w$ containing a vertex which belongs to $B$.

Suppose that $Q$ has length less than $d_r$. Then $Q$ is an $r$-string and the two endpoints of $Q$ do not belong to the same block of $\pi$ (else $\Gamma_w/\pi$ contains an $r$-loop of length less than $d_r$, contradicting the assumption that it is strongly admissible). In particular, at least one of the two endpoints of $Q$ does not belong to $B$, so the block of $\pi$ containing this endpoint belongs to no path of $\Gamma_w/\pi$ other than $P$. This implies that this endpoint belongs only to paths of $\Gamma_w$ which have color $r$, which is plainly impossible by the definition of $\Gamma_w$.

Thus $Q$ has length at least $d_r$, from which it is clear that $d_r<\infty$ and some cyclic permutation of the symbols of $w$ is of the form $vu_r^{d_r}$ or $v\left(u_r^*\right)^{d_r}$.
\end{proof}
\end{lemma}
\begin{prop}\label{prop scon chi word-graphs}
If \F\ is as in \ref{notn s}(b) and $w\in\F\setminus\{e\}$, then
\begin{equation}
\mathrm{card}\big\{\pi\in\SCon(\Gamma_w)\,\big|\,\chi(\Gamma_w/\pi)=1\big\}=
\begin{cases}1&\text{if $w\approx e$}
\\0&\text{if $w\not\approx e$.}
\end{cases}
\end{equation}
\begin{proof}
The proof is by induction on the length of $w$. If $w$ has length 1, then $w=g_r$ or $g_r^*$ for some $r\in[s]$, $\Gamma_w$ is an $r$-loop of length 1, and $\pi:=\big\{\{1\}\big\}$ is the only partition of the vertex set of $\Gamma_w$. It is clear that $\pi\in\Con(\Gamma_w)$, that $\chi(\Gamma_w/\pi)=1$, and that $\pi\in\SCon(\Gamma_w)$ iff $d_r=1$, and hence iff $w\approx e$.

Now suppose that the assertion holds for every word in $\F\setminus\{e\}$ having length less than or equal to some positive integer $n$, and let $w\in\F\setminus\{e\}$ have length $n+1$.

First suppose that there exists $\pi\in\SCon(\Gamma_w)$ such that $\chi(\Gamma_w/\pi)=1$. By Lemma \ref{lemma paths}, there exists a path $P$ of $\Gamma_w/\pi$ such that $(\Gamma_w/\pi)\setminus P$ has at most one non-trivial component. On the other hand, since $\chi(\Gamma_w/\pi)=1$, it follows from Proposition \ref{connected chi} that $P$ has exactly one vertex in common with each connected component of $(\Gamma_w/\pi)\setminus P$. Thus at most one vertex of $P$ belongs to other paths of $\Gamma_w/\pi$, so by Lemma \ref{lemma cyclic} some cyclic permutation of the symbols of $w$ is of the form $vx$ with $x$ as in \eqref{x}. To show that $w\approx e$, it suffices to show that $vx\approx e$. If $v=e$, this is immediate. If $v\neq e$, then there exists $\pi\in\SCon(\Gamma_{vx})$ such that $\chi(\Gamma_{vx}/\pi)=1$ (since $\Gamma_{vx}$ is isomorphic to $\Gamma_w$), so by Lemma \ref{lemma scon chi word-graphs}(b) there exists
$\pi\in\SCon(\Gamma_v)$ such that $\chi(\Gamma_v/\pi)=1$; thus $v\approx e$ by the induction hypothesis, so $vx\approx e$.

Next suppose that $w\approx e$. Then some cyclic permutation of the symbols of $w$ is of the form $vx$ with $v\approx e$ and $x$ as in \eqref{x}. Since $\Gamma_{vx}$ is isomorphic to $\Gamma_w$, it suffices to show that
$$\card\big\{\pi\in\SCon(\Gamma_{vx})\,\big|\,\chi(\Gamma_{vx}/\pi)=1\big\}=1.$$
If $v=e$, this follows from Lemma \ref{lemma scon chi word-graphs}(a), whereas if $v\neq e$, it follows from Lemma \ref{lemma scon chi word-graphs}(b), since $\card\big\{\rho\in\SCon(\Gamma_v)\,\big|\,\chi(\Gamma_v/\rho)=1\big\}=1$ by the
induction hypothesis.
\end{proof}
\end{prop}
\section{Permutations with Restricted Cycle Lengths}\label{sec random
permutations}
\begin{permutations}\label{permutations}
Let $A$ be a non-empty set of positive integers. For every positive integer $N$,
$\S_N$ denotes the set of all permutations of $[N]$ and $\S_N\toA$ denotes the
subset of $\S_N$ defined by
\begin{equation}
\S_N\toA=\{\sigma\in\S_N\,|\,\text{the length of every cycle of $\sigma$
belongs to $A$}\}.
\end{equation}
For every positive integer $N$, the cardinality of the set $\S_N\toA$ is denoted by $a_N\toA$; by convention, $a_0\toA:=1$. It is well known (see, for instance, \cite[Theorem 3.53]{Bona}) that the exponential generating function for the
numbers $a_N\toA$ satisfies
\begin{equation}\label{egf}
\sum_{N=0}^\infty\frac{a_N\toA}{N!}z^N=\exp\left(\sum_{k\in A}\frac{z^k}{k}\right).
\end{equation}
It is clear that $a_N\toA=0$ for every $N$ which is not a multiple of $\gcd(A)$, where $\gcd(A)$ denotes the
greatest common divisor of the elements of $A$. On the other hand, we have the following.
\end{permutations}
\begin{lemma}\label{lemma gcd}
If $A$ is a non-empty set of positive integers, then $a_N\toA>0$ for every sufficiently large multiple $N$ of $\gcd(A)$.
\begin{proof}
Let $D=\gcd(A)$. Then there exist $k_1,k_2,\ldots,k_n\in A$ and integers $\alpha_1,\alpha_2,\ldots,\alpha_n$ such that
$$D=\sum_{i=1}^n\alpha_ik_i.$$
Let $N=\frac{k_1}{D}\sum\limits_{i=1}^N\left|\alpha_i\right|k_i$. For every positive integer $m$, let $q,r$ be non-negative integers such that $m=q\left(\frac{k_1}{D}\right)+r$ and $0\leq r<\frac{k_1}{D}$; then
$$N+mD=N+qk_1+rD=qk_1+\sum_{i=1}^n\left(\frac{k_1}{D}\left|\alpha_i\right|+r\alpha_i\right)k_i.$$
Thus every multiple of $D$ which is greater than $N$ can be written as a linear combination of elements of $A$ with positive coefficients, from which the lemma follows.
\end{proof}
\end{lemma}
\begin{cycles}\label{cycles}
Let $A$ be a non-empty set of positive integers, let $k\in A$, and let $N$ be a positive integer such that $N\geq k$ and $\S_N\toA\neq\emptyset$. For every $i,j\in[N]$, it is easily seen that the map
\begin{align*}
\left.\left\{\sigma\in\S_N\toA\right|\text{$i$ belongs to a cycle of length $k$}\right\}&\longrightarrow\left.\left\{\sigma\in\S_N\toA\right|\text{$j$ belongs to a cycle of length $k$}\right\}
\\\sigma&\longmapsto(i,j)\sigma(i,j)
\end{align*}
is a bijection, so the quantity $p_N\toA(k)$ defined by
\begin{equation}
p_N\toA(k)=\frac{\card\left.\left\{\sigma\in\S_N\toA\,\right|\,\text{$i$ belongs to a cycle of length $k$}\right\}}{a_N\toA}
\end{equation}
(notation as in \ref{permutations}) is independent of the choice of $i\in[N]$. Clearly, $p_N\toA(k)$ is the probability that, in a randomly chosen permutation from $\S_N\toA$, $i$ belongs to a cycle of length $k$; equivalently, $p_N\toA(k)$ is $\frac{k}{N}$ times the expected number of cycles of length $k$ of a randomly chosen permutation from $\S_N\toA$. Note that for every $i\in[N]$, we have that
\begin{equation}\label{cardsnai}
\card\big\{\sigma\in\S_N\toA\,\big|\,\text{$i$ belongs to a cycle of length $k$}\big\}=(N-1)(N-2)\cdots(N-k+1)a_{N-k}\toA,
\end{equation}
because there are $N-1$ equally likely possibilities for $\sigma(i)$, and for each of these possibilities there are $N-2$ equally likely possibilities for $\sigma^2(i)$ etc., and once the cycle (of length $k$) to which $i$ belongs is determined, there are $a_{N-k}\toA$ possibilities for the remaining cycles. Thus
\begin{equation}\label{pnk}
p_N\toA(k)=\frac{(N-1)!\;a_{N-k}\toA}{(N-k)!\;a_N\toA}=\frac{1}{N}\cdot\frac{a_{N-k}\toA/(N-k)!}{a_N\toA/N!}.
\end{equation}
\end{cycles}
\begin{remks}
Let $A$ be a non-empty set of positive integers and let $k\in A$. We shall be interested in the asymptotic behavior of $p_N\toA(k)$ as $N\to\infty$ with $\S_N\toA\neq\emptyset$. This asymptotic behavior is described in Proposition \ref{limpnk} below, under the assumption that the set $A$ satisfies the following condition:
\begin{equation}\label{condition}
\text{Either $A$ is finite, or $\sum_{j\in\N\setminus A}\frac{1}{j}<\infty$.}
\end{equation}
[Note that the second hypothesis on $A$ implies that $\gcd(A)=1$, so $S_N\toA\neq\emptyset$ for all sufficiently large $N$.] The part of Proposition \ref{limpnk} which deals with the case when $A$ is finite relies on Lemma \ref{Hayman} below, which uses Hayman's method for determining the asymptotics of power series coefficients of certain complex analytic functions. In turn, Lemma \ref{Hayman} relies on the following intermediate step.
\end{remks}
\begin{lemma}\label{analytic facts}
Let $\mathcal N$ be a neighborhood of the origin in the complex plane, let $F:\mathcal N\to\C$ be an analytic function with $F(0)=1$, let $k$ be a positive integer, let $C>0$, and for every positive integer $N$ with $t_N:=\left(\frac{C}{N}\right)^{1/k}\in\mathcal N$, let $r_N=\frac{F(t_N)}{t_N}$. Then:

(a) $\lim\limits_{N\to\infty}\left(r_N^k-r_{N-1}^k\right)=\frac{1}{C}.$

(b) $\lim\limits_{N\to\infty}\left(r_N^j-r_{N-1}^j\right)=0$ for every $0\leq j<k.$

(c) $\frac{r_N^N}{r_{N-1}^{N-1}}\sim\left(\frac{\mathrm eN}{C}\right)^{1/k}$ as $N\to\infty.$
\begin{proof}
The proof uses the simple observation that
$$\lim_{N\to\infty}\bigl(N^\alpha-(N-1)^\alpha\bigr)=0\qquad\text{for every }0\leq\alpha<1,$$
which immediately implies that
\begin{equation}\label{observation}
\lim_{N\to\infty}\left(\frac{1}{t_N^m}-\frac{1}{t_{N-1}^m}\right)=0\qquad\text{for every }0\leq m<k.
\end{equation}
To prove (a) and (b), let $0\leq j\leq k$ and consider the power series expansion of $F^j$
around the origin:
$$F(t)^j=:1+\sum\limits_{i=1}^j\alpha_it^i+t^{j+1}G_j(t).$$
Then for every sufficiently large $N$, we have that
$$r_N^j=\frac{F(t_N)^j}{t_N^j}=\frac{1}{t_N^j}+\sum_{i=1}^j\frac{\alpha_i}{t_N^{j-i}}+t_NG_j(t_N),$$
so \eqref{observation} implies that
$$\lim_{N\to\infty}\left(r_N^j-r_{N-1}^j\right)=\lim_{N\to\infty}\left(\frac{1}{t_N^j}-\frac{1}{t_{N-1}^j}\right)=\begin{cases}\frac{1}{C}&\text{if $j=k$}\\0&\text{if $j<k$.}
\end{cases}$$
To prove (c), consider the power series expansion of $\log F$ around the origin:
$$\log F(t)=:\sum\limits_{i=1}^k\beta_it^i+t^{k+1}G(t).$$
Then for every sufficiently large $N$, we have that
$$F(t_N)^N=\exp\big(N\log F(t_N)\big)=\exp\left(\frac{C}{t_N^k}\log F(t_N)\right)=\exp\left(\sum_{i=1}^k\frac{C\beta_i}{t_N^{k-i}}+Ct_NG(t_N)\right),$$
so \eqref{observation} implies that
$$\lim_{N\to\infty}\frac{F(t_N)^N}{F(t_{N-1})^{N-1}}=1,$$
from which it follows that
$$\frac{r_N^N}{r_{N-1}^{N-1}}\sim\frac{t_{N-1}^{N-1}}{t_N^N}=\left(\frac{N^N}{C(N-1)^{N-1}}\right)^{1/k}\sim\left(\frac{\mathrm eN}{C}\right)^{1/k}\qquad\text{as $N\to\infty$.}$$
\end{proof}
\end{lemma}
\begin{lemma}\label{Hayman}
Let $n$ be a positive integer, let $c_1,c_2,\ldots,c_n>0$, let $(k_i)_{i=1}^n$ be an increasing sequence of positive integers, let $\sum\limits_{N=0}^\infty b_Nw^N$ be the power series expansion around the origin of the function $f(w):=\exp\left(\sum\limits_{i=1}^nc_iw^{k_i}\right),$ and suppose that $b_N>0$ for all sufficiently large $N$. Then
\begin{equation}\label{eqHayman}
\frac{b_{N-1}}{b_N}\sim\left(\frac{N}{c_nk_n}\right)^{1/k_n}\qquad\text{as $N\to\infty$.}
\end{equation}
\begin{proof}
The proof relies on Hayman's method for determining the asymptotics of power series coefficients of certain complex analytic functions, which is described in detail in \cite[Section 5.4]{Wilf}. By Theorem 5.4.1 of that work, we have that
\begin{equation}\label{Haymaneq}
b_N\sim\frac{f(r_N)}{r_N^N\sqrt{2\pi\left(\sum\limits_{i=1}^nc_ik_i^2r_N^{k_i}\right)}}\qquad\text{as $N\to\infty$,}
\end{equation}
where $r_N$ is the positive real root of the equation
\begin{equation}\label{eq1}
\sum\limits_{i=1}^nc_ik_ir_N^{k_i}=N.
\end{equation}
It is immediate that $\lim\limits_{N\to\infty}r_N=\infty$, and hence that
$$\sum\limits_{i=1}^nc_ik_i^2r_N^{k_i}\sim c_nk_n^2r_N^{k_n}\sim k_n\left(\sum\limits_{i=1}^nc_ik_ir_N^{k_i}\right)=k_nN\qquad\text{as $N\to\infty$},$$
so \eqref{Haymaneq} implies that
\begin{equation}\label{eq2}
\frac{b_{N-1}}{b_N}\sim\frac{f(r_{N-1})}{f(r_N)}\cdot\frac{r_N^N}{r_{N-1}^{N-1}}=\frac{r_N^N}{r_{N-1}^{N-1}}\exp\left(-\sum_{i=1}^nc_i\left(r_N^{k_i}-r_{N-1}^{k_i}\right)\right)\qquad\text{as $N\to\infty$}.
\end{equation}
To evaluate the right-hand side of \eqref{eq2}, further analysis of \eqref{eq1} is needed. First note that \eqref{eq1} can be put in the form
$$\frac{1}{r_N}=\left(\frac{c_nk_n}{N}\right)^{1/k_n}\left(1+\sum_{i=1}^{n-1}\frac{c_ik_i}{c_nk_nr_N^{k_n-k_i}}\right)^{1/k_n}.$$
Via the generalized binomial expansion of the right-hand side and with the substitutions $u_N:=\frac{1}{r_N}$ and $t_N:=\left(\frac{c_nk_n}{N}\right)^{1/k_n}$, this can in turn be written as
$$u_N=t_N\left(1+\sum_{j=1}^\infty\alpha_ju_N^j\right)$$
for some coefficients $\alpha_j$ such that the series $\varphi(u):=1+\sum\limits_{j=1}^\infty\alpha_ju^j$ has positive radius of convergence. Since the function $\frac{u}{\varphi(u)}$ is conformal for sufficiently small $u$, it has an analytic inverse $g$ [with $g(0)=0$ and $g'(0)=1$] in a small neighborhood of the origin. Thus $u_N=g(t_N)$ for all sufficiently large $N$. Letting $F(t):=\frac{t}{g(t)}$ for all sufficiently small $t$ and $F(0):=1$, it follows that $F$ is analytic and $r_N=\frac{F(t_N)}{t_N}$ for all large enough $N$,
so \eqref{eq2} and Lemma \ref{analytic facts} imply that
$$\frac{b_{N-1}}{b_N}\sim\left(\frac{\mathrm eN}{c_nk_n}\right)^{1/k_n}\mathrm\exp(-c_n/c_nk_n)=\left(\frac{N}{c_nk_n}\right)^{1/k_n}\qquad\text{as $N\to\infty$.}$$
\end{proof}
\end{lemma}
\begin{prop}\label{limpnk}
Let $A$ be a non-empty set of positive integers which satisfies condition \eqref{condition}, let $k\in A$, let $d=\sup A$, and for every positive integer $N$ with $S_N\toA\neq\emptyset$, let $p_N\toA(k)$ be as in \ref{cycles}. Then
\begin{equation}
p_N\toA(k)\sim N^{(k/d)-1}\qquad\text{as $N\to\infty$ with $S_N\toA\neq\emptyset$}
\end{equation}
(with the convention that $k/\infty:=0$).
\begin{proof}
Suppose first that $A$ is finite and that the elements of $A$ are $k_1,k_2,\ldots,k_n$, written in increasing order (so $d=k_n$). Let $D=\gcd(A)$. Then
\begin{equation*}
\sum_{N=0}^\infty\frac{a_{DN}\toA}{(DN)!}z^{DN}=\exp\left(\sum_{i=1}^n\frac{z^{k_i}}{k_i}\right),\tag*{[by \eqref{egf}]}\end{equation*}
which, via the substitution $w:=z^D$, can be rewritten as
$$\sum_{N=0}^\infty\frac{a_{DN}\toA}{(DN)!}w^N=\exp\left(\sum_{i=1}^n\frac{w^{k_i/D}}{k_i}\right).$$
Since $a_{DN}\toA>0$ for all sufficiently large $N$ by Lemma \ref{lemma gcd}, it follows from Lemma \ref{Hayman} that
$$\frac{a_{DN-D}\toA/(DN-D)!}{a_{DN}\toA/(DN)!}\sim(DN)^{D/k_n}\qquad\text{as $N\to\infty$.}$$
Hence \eqref{pnk} implies that
$$p_{DN}\toA(k)=\frac{1}{DN}\cdot\frac{a_{DN-k}\toA/(DN-k)!}{a_{DN}\toA/(DN)!}\sim(DN)^{(k/k_n)-1}\qquad\text{as $N\to\infty$.}$$
Next suppose that $\sum\limits_{j\in\N\setminus A}\frac{1}{j}<\infty$ [which implies that $d=\infty$ and that $\gcd(A)=1$, so $S_N\toA\neq\emptyset$ for all sufficiently large $N$]. For every $N\geq1$ we have that
\begin{equation*}
a_N\toA=\sum_{j\in A\cap[N]}\card\big\{\sigma\in\S_N\toA\,\big|\,\text{1 belongs to a cycle of length $j$}\big\}=\sum_{j\in
A\cap[N]}\frac{(N-1)!\;a_{N-j}\toA}{(N-j)!}\tag*{\text{[by \eqref{cardsnai}],}}
\end{equation*}
so
$$\frac{a_N\toA}{N!}=\frac{1}{N}\sum_{j\in A\cap[N]}\frac{a_{N-j}\toA}{(N-j)!}.$$
Because of the assumption $\sum\limits_{j\in\N\setminus A}\frac{1}{j}<\infty$, it now follows immediately from Theorem 1 of \cite{Hildebrand} that
$$\lim_{N\to\infty}\frac{a_N\toA}{N!}=\exp\left(-\sum_{j\in\N\setminus A}\frac{1}{j}\right).$$
(This conclusion also follows from Theorem 1.1 of \cite{Manstavicius}, in the more general context of combinatorial assemblies.) Hence \eqref{pnk} implies that
$$p_N\toA(k)=\frac{1}{N}\cdot\frac{a_{N-k}\toA/(N-k)!}{a_N\toA/N!}\sim\frac{1}{N}\qquad\text{as $N\to\infty$.}$$
\end{proof}
\end{prop}
\begin{remks}
(a) Condition \eqref{condition} is not necessary for the conclusion of Proposition \ref{limpnk} to hold, as shown by the following example.

Let $D$ be a positive integer and let $A=\{DN\,|\,N\geq1\}$. By \eqref{egf}, we have that
\begin{align*}\sum_{N=0}^\infty\frac{a_{DN}\toA}{(DN)!}z^{DN}&=\exp\left(\sum_{N=1}^\infty\frac{z^{DN}}{DN}\right)=\exp\left(-\frac{1}{D}\log\big(1-z^D\big)\right)=\left(1-z^D\right)^{-1/D}
\\&=\sum_{N=0}^\infty\binom{-1/D}{N}\big(-z^D\big)^N,
\end{align*}
where $\binom{-1/D}{N}:=\frac{1}{N!}\left(-\frac{1}{D}\right)\left(-\frac{1}{D}-1\right)\cdots\left(-\frac{1}{D}-N+1\right)$ is a generalized binomial coefficient. Thus
$$\frac{a_{DN}\toA}{(DN)!}=\frac{1}{N!}\left(\frac{1}{D}\right)\left(\frac{1}{D}+1\right)\cdots\left(\frac{1}{D}+N-1\right)\qquad\forall N\geq1,$$
from which it immediately follows that
$$\lim_{N\to\infty}\frac{a_{DN-D}\toA/(DN-D)!}{a_{DN}\toA/(DN)!}=1.$$
Hence \eqref{pnk} implies that for every $k\in A$, we have that
$$p_{DN}\toA(k)=\frac{1}{DN}\cdot\frac{a_{DN-k}\toA/(DN-k)!}{a_{DN}\toA/(DN)!}\sim\frac{1}{DN}\qquad\text{as $N\to\infty$.}$$

(b) On the other hand, the following example shows the existence of sets of positive integers for which the conclusion of Proposition \ref{limpnk} does not hold.

First let us construct inductively a sequence $(N_k)_{k=1}^\infty$ of positive integers such that
$$N_{k+1}p_{N_{k+1}}^{(A_k)}(1)>2\qquad\forall k\geq1,$$
where $A_k:=\{1,N_1+1,N_2+1,\ldots,N_k+1\}$. Let $N_1$ be any positive integer. If $k\geq1$ and $N_1,N_2,\ldots,N_k$ have been chosen, Proposition \ref{limpnk} implies that $p_N^{(A_k)}(1)\sim N^{1/(N_k+1)-1}$ as $N\to\infty$. In particular, $Np_N^{(A_k)}(1)\to\infty$ as $N\to\infty$, so there exists $N_{k+1}>N_k$ such that $N_{k+1}p_{N_{k+1}}^{(A_k)}>2$.

Now let $A=\bigcup\limits_{k=1}^\infty A_k$. Since $p_N\toA(1)$ depends only on the elements of $A$ which are less than or equal to $N$, we have that $p_{N_{k+1}}\toA(1)=p_{N_{k+1}}^{(A_k)}(1)\;\forall k\geq1$. Thus $N_{k+1}p_{N_{k+1}}\toA(1)>2$ $\forall k\geq1$, so $p_N\toA(1)\not\sim\frac{1}{N}$ as $N\to\infty$.
\end{remks}
\begin{permutations and graphs}\label{permutations and graphs}
(a) For every non-empty set $A$ of positive integers, an \emph{$A$-graph} is a monocolored admissible graph $\Gamma$ with the properties that $l(L)\in A$ for every loop $L$ of $\Gamma$ and $l(S)<\sup A$ for every string $S$ of $\Gamma$ (notation and terminology as in \ref{admissible graphs} and \ref{loop-characteristic}).

(b) Let $A$ be a non-empty set of positive integers, let $N$ be a positive integer such that $\S_N\toA\neq\emptyset$, and let $\Gamma$ be a monocolored graph with vertex set contained in $[N]$. A permutation $\sigma\in\S_N\toA$ is said to be \emph{compatible with $\Gamma$} iff for every edge $i\rightarrow j$ of $\Gamma$, we have that $\sigma(i)=j$. The set of all permutations in $\S_N\toA$ which are compatible with $\Gamma$ is denoted by $\S_N\toA(\Gamma)$. It is clear that, if $\Gamma$ is not an $A$-graph, then $\S_N\toA(\Gamma)=\emptyset$.

(c) Let $A$ be a non-empty set of positive integers, let $N$ be a positive integer such that $\S_N\toA\neq\emptyset$, and let $\Gamma$ be an $A$-graph with vertex set $V$ consisting of at most $N$ vertices $v_1,v_2,\ldots,v_k$. For every injective function $f:V\to[N]$, let $f(\Gamma)$ denote the graph obtained from $\Gamma$ by relabeling $v_1,v_2,\ldots,v_k$ as $f(v_1),f(v_2),\ldots,f(v_k)$, respectively. If $f$ and $g$ are two such functions, it is easily seen that the map
\begin{align*}
\S_N\toA\big(f(\Gamma)\big)&\longrightarrow\S_N\toA\big(g(\Gamma)\big),
\\\sigma&\longmapsto\big(f(v_1),g(v_1)\big)\cdots\big(f(v_k),g(v_k)\big)\sigma\big(f(v_1),g(v_1)\big)\cdots\big(f(v_k),g(v_k)\big)
\end{align*}
is a bijection, so the quantity $p_N\toA(\Gamma)$ defined by
\begin{equation}
p_N\toA(\Gamma)=\frac{\card\;\S_N\toA\big(f(\Gamma)\big)}{a_N\toA}
\end{equation}
is independent of the choice of the injective function $f:V\to[N]$. Clearly, $p_N\toA(\Gamma)$ is the probability that a randomly chosen permutation from $\S_N\toA$ is compatible with $f(\Gamma)$.
\end{permutations and graphs}
\begin{prop}\label{limpng}
Let $A$ be a non-empty set of positive integers which satisfies condition \eqref{condition}, let $d=\sup A$, let $\Gamma$ be an $A$-graph [terminology as in \ref{permutations and graphs}(a)], and for every positive integer $N$ with $S_N\toA\neq\emptyset$, let $p_N\toA(\Gamma)$ be as in \ref{permutations and graphs}(c). Then
\begin{equation}
p_N\toA(\Gamma)\sim N^{-E(\Gamma)+L(\Gamma)+\sum\limits_{\substack{\text{loops}\\\text{$P$ of $\Gamma$}}}\left[\frac{l(P)}{d}-1\right]}\qquad\text{as $N\to\infty$ with $\S_N\toA\neq\emptyset$}
\end{equation}
\Big(with the convention that $\frac{l(P)}{\infty}:=0$\Big).
\begin{proof}
The proof  is by induction on the number of paths of $\Gamma$. If $\Gamma$ has no paths (and therefore no edges), then it easily seen that $p_N\toA(\Gamma)=1$ for every positive integer $N$ with $\S_N\toA\neq\emptyset$.

Now suppose that the assertion of the proposition holds for every $A$-graph with number of paths less than some positive integer $n$, and let $\Gamma$ be an $A$-graph with $n$ paths. Without loss of generality, we may assume that the vertex set of $\Gamma$ is contained in $[N_0]$ for some positive integer $N_0$. Let $Q$ be any path of $\Gamma$ and let $\Gamma\setminus Q$ be the graph obtained from $\Gamma$ by removing all edges of $Q$.

First suppose that $Q$ is a loop $v_1\rightarrow v_2\rightarrow\cdots\rightarrow v_k\rightarrow v_1$ (of some length $k\in A$). Then for every $N\geq N_0$ with $\S_N\toA\neq\emptyset$, we have that
$$p_N\toA(\Gamma)=p_N\toA(Q)\cdot p_{N-k}\toA(\Gamma\setminus Q)$$
and that
$$p_N\toA(Q)=\frac{p_N\toA(k)}{(N-1)(N-2)\cdots(N-k+1)},$$
because, given that $v_1$ belongs to a cycle of length $k$ of a permutation $\sigma\in\S_N\toA$, there are $N-1$ equally likely possibilities for $\sigma(v_1)$, and for each of these possibilities there are $N-2$ equally likely possibilities for $\sigma^2(v_1)$ etc. Then Proposition \ref{limpnk} and the induction hypothesis imply that
\begin{align*}
p_N\toA(\Gamma)&\sim N^{-k+1}N^{(k/d)-1}N^{-E(\Gamma\setminus Q)+L(\Gamma\setminus Q)+\sum\limits_{\substack{\text{loops $P$}
\\\text{of $\Gamma\setminus Q$}}}\left[\frac{l(P)}{d}-1\right]}
\\&=N^{-E(\Gamma)+L(\Gamma)+\sum\limits_{\substack{\text{loops}
\\\text{$P$ of $\Gamma$}}}\left[\frac{l(P)}{d}-1\right]}\qquad\text{as $N\to\infty$ with $S_N\toA\neq\emptyset$}.
\end{align*}
Finally, suppose that $Q$ is a string $v_1\rightarrow v_2\rightarrow\cdots\rightarrow v_k\rightarrow v_{k+1}$ (of some length $k<\sup A$). Let $V$ be the set of vertices which belong to the paths of $\Gamma$ other than $Q$, and let $|V|$ be the number of vertices in $V$. Let $q_N$ be the probability that a randomly chosen permutation $\sigma\in\S_N\toA(\Gamma\setminus Q)$ has the property that the numbers $\sigma(v_1),\sigma^2(v_1),\ldots,\sigma^k(v_1)$ are distinct from each other and from the vertices in $V$. Then for every $N\geq N_0$ with $\S_N\toA\neq\emptyset$, we have that
$$p_N\toA(\Gamma)=\frac{p_N\toA(\Gamma\setminus Q)\cdot q_N}{(N-|V|-1)(N-|V|-2)\cdots(N-|V|-k)},$$
because there are $N-|V|-1$ equally likely possibilities for $\sigma(v_1)$, and for each of these possibilities there are $N-|V|-2$ equally likely possibilities for $\sigma^2(v_1)$ etc. Then the induction hypothesis implies that
\begin{align*}
p_N\toA(\Gamma)&\sim q_N\cdot N^{-k}N^{-E(\Gamma\setminus Q)+L(\Gamma\setminus Q)+\sum\limits_{\substack{\text{loops $P$}
\\\text{of $\Gamma\setminus Q$}}}\left[\frac{l(P)}{d}-1\right]}
\\&=q_N\cdot N^{-E(\Gamma)+L(\Gamma)+\sum\limits_{\substack{\text{loops}
\\\text{$P$ of $\Gamma$}}}\left[\frac{l(P)}{d}-1\right]}\qquad\text{as $N\to\infty$ with $S_N\toA\neq\emptyset$},
\end{align*}
so it suffices to show that $q_N\to1$ (or equivalently, that $1-q_N\to0$) as $N\to\infty$ with $\S_N\toA\neq\emptyset$.

Let $N\geq N_0$ be such that $\S_N\toA\neq\emptyset$. Given a permutation $\sigma\in\S_N\toA(\Gamma\setminus Q)$, there are only two (mutually exclusive) ways in which $\sigma$ can fail to have the property that the numbers $\sigma(v_1),\sigma^2(v_1),\ldots,\sigma^k(v_1)$ are distinct from each other and from the vertices in $V$; namely, either
\pagebreak

(1) at least one of $\sigma(v_1),\sigma^2(v_1),\ldots,\sigma^k(v_1)$ belongs to $V$, or
\vspace{2mm}

(2) $v_1$ belongs to a cycle of $\sigma$ of length at most $k$ which does not contain any vertex in $V$.
\vspace{2mm}

For any vertex $v\in V$ and for every $i\in[k]$, the number of possibilities for $\sigma^{-i}(v)$ is either 1 or $N-|V|$ (and these possibilities are equally likely), so the probability $q_N^{(1)}$ of (1) occurring satisfies
$$q_N^{(1)}\leq\sum\limits_{v\in V}\sum\limits_{i=1}^k\frac{1}{N-|V|}=\frac{k\cdot|V|}{N-|V|}.$$
On the other hand, the probability $q_N^{(2)}$ of (2) occurring satisfies
$$q_N^{(2)}\leq\sum\limits_{j\in A\cap[k]}\frac{p_N\toA(j)\cdot p_{N-j}\toA(\Gamma\setminus Q)}{p_N\toA(\Gamma\setminus Q)},$$
which, by the induction hypothesis and Proposition \ref{limpnk} (together with the fact that $k<\sup A$), converges to 0 as $N\to\infty$ with $\S_N\toA\neq\emptyset$. 

Thus $1-q_N=q_N^{(1)}+q_N^{(2)}\to0$ as $N\to\infty$ with $\S_N\toA\neq\emptyset.$
\end{proof}
\end{prop}
\section{Asymptotic Freeness.}
\begin{non-commutative probability spaces}\label{non-commutative probability spaces}
(a) A \emph{non-commutative $*$-probability space} (or, simply, a \emph{$*$-probability space}) is a pair $(\A,\varphi)$, where \A\ is a unital $*$-algebra (the unit of which is denoted by 1) and $\varphi:\A\to\C$ is a linear functional which satisfies
\begin{equation}
\varphi(1)=1\qquad\text{and}\qquad\varphi(a^*a)\geq0\quad\forall a\in\A.
\end{equation}
Elements of \A\ are called \emph{non-commutative random variables} (or, simply, \emph{random variables}).

(b) Let $(\A,\varphi)$ be a $*$-probability space. A family $(\A_j)_{j\in J}$ of unital subalgebras of \A\ is said to be \emph{free} iff for every integer $n\geq1$, for every $j_1,j_2,\ldots,j_n\in J$ such that $j_1\neq j_2\neq\cdots\neq j_n$, and for every $a_1\in\A_{j_1}$, $a_2\in\A_{j_2}$,\ldots,$a_n\in\A_{j_n}$ such that $\varphi(a_1)=\varphi(a_2)=\cdots=\varphi(a_n)=0$, we have that $\varphi(a_1a_2\cdots a_n)=0$. More generally, a family $(a_j)_{j\in J}$ of random variables in \A\ is said to be \emph{free} (respectively, \emph{$*$-free}) iff the family $(\A_j)_{j\in J}$ is free, where each $\A_j$ is the unital subalgebra (respectively, $*$-subalgebra) of \A\ generated by $a_j$.

(c) Let $(\A,\varphi)$ be a $*$-probability space. If $d$ is a positive integer, a \emph{$d$-Haar unitary} in $(\A,\varphi)$ is a unitary random variable $u\in\A$ such that $u^d=1$ and $\varphi(u^n)=0$ whenever $n$ is not a multiple of $d$. An \emph{$\infty$-Haar unitary} in $(\A,\varphi)$ is a unitary random variable $u\in\A$ such that $\varphi(u^n)=0\;\forall n\neq0$.
\end{non-commutative probability spaces}
\begin{random matrices}\label{random matrices}
(a) For every positive integer $N$, the set of all $N\times N$ matrices with complex entries is denoted by $\M _N$. If $A\in\M_N$, then $\tr\;A$ denotes the \emph{normalized trace} of $A$:
\begin{equation}
\tr\;A:=\frac{1}{N}\sum_{i=1}^N A_{ii}.
\end{equation}

(b) We shall work with a fixed probability space $(\Omega,\mathcal F,P)$, over which we shall consider \emph{random variables} and \emph{random $N\times N$ matrices}; that is, measurable functions $f:\Omega\to\C$ and $A:\Omega\to\M_N$. The \emph{expectation} of an integrable random variable $f$ over $(\Omega,\mathcal F,P)$ is
\begin{equation}
\E(f):=\int_\Omega f(\omega)\dif P(\omega).
\end{equation}

(c) Let $(a_r)_{r=1}^s$ be a family of random variables in a $*$-probability space $(\A,\varphi)$, let $(N_k)_{k=1}^\infty$ be an increasing sequence of positive integers, and for each $k$ let $\left(A^{(k)}_r\right)_{r=1}^s$ be a family of random $N_k\times N_k$ matrices over $(\Omega,\mathcal F,P)$. The families $\left(A^{(k)}_r\right)_{r=1}^s$ are said to \emph{converge in $*$-distribution} to $(a_r)_{r=1}^s$ as $k\to\infty$, written $\left(A^{(k)}_r\right)_{r=1}^s\converges(a_r)_{r=1}^s$, iff 
\begin{equation}
\lim_{k\to\infty}\E\left[\tr\;\left(A^{(k)}_{r_1}\right)^{\varepsilon_1}\left(A^{(k)}_{r_2}\right)^{\varepsilon_2}\cdots\left(A^{(k)}_{r_n}\right)^{\varepsilon_n}\right]=\varphi\left(a_{r_1}^{\varepsilon_1}a_{r_2}^{\varepsilon_2}\cdots a_{r_n}^{\varepsilon_n}\right)
\end{equation}
for every positive integer $n$, for every $r_1,r_2,\ldots,r_n\in[s]$, and for every $\varepsilon_1,\varepsilon_2,\ldots,\varepsilon_n\in\{1,*\}$.

The families $\left(A^{(k)}_r\right)_{r=1}^s$ are said to \emph{converge in $*$-distribution almost surely} to $(a_r)_{r=1}^s$ as $k\to\infty$, written $\left(A^{(k)}_r\right)_{r=1}^s\convergesas(a_r)_{r=1}^s$, iff they satisfy the
stronger condition that
\begin{equation}
\lim_{k\to\infty}\tr\;\left(A^{(k)}_{r_1}\right)^{\varepsilon_1}\left(A^{(k)}_{r_2}\right)^{\varepsilon_2}\cdots\left(A^{(k)}_{r_n}\right)^{\varepsilon_n}=\varphi\left(a_{r_1}^{\varepsilon_1}a_{r_2}^{\varepsilon_2}\cdots a_{r_n}^{\varepsilon_n}\right)\qquad\text{almost surely}
\end{equation}
for every positive integer $n$, for every $r_1,r_2,\ldots,r_n\in[s]$, and for every $\varepsilon_1,\varepsilon_2,\ldots,\varepsilon_n\in\{1,*\}$.
\end{random matrices}
In many concrete instances where convergence in $*$-distribution is known to hold, it can be shown that it holds almost surely by means of the following fact, the proof of which can be found, for instance, embedded in the proof of Corollary 3.9 in \cite{Thorbjornsen}.
\begin{prop}\label{prop as}
Let $\lambda\in\C$ and let $\left(f_k\right)_{k=1}^\infty$ be a sequence of random variables over $(\Omega,\mathcal F,P)$ such that $\lim\limits_{k\to\infty}\E(f_k)=\lambda$ and such that
\begin{equation}
\sum_{k=1}^\infty\left[\E\big(\left|f_k\right|^2\big)-\big|\E(f_k)\big|^2\right]<\infty.
\end{equation}
Then $\lim\limits_{k\to\infty}f_k=\lambda$ almost surely.
\end{prop}
\begin{random permutations}\label{random permutations}
(a) For every positive integer $N$ and for every $\sigma\in\S_N$, $M_\sigma$ denotes the $N\times N$ matrix defined by
\begin{equation}
(M_\sigma)_{i,j}=
\begin{cases}
1&\text{if $\sigma(j)=i$}
\\0&\text{if $\sigma(j)\neq i$}	     
\end{cases}\qquad\forall i,j\in[N].
\end{equation}

(b) Let $N$ be a positive integer. A \emph{random $N\times N$ permutation matrix} over \OFP\ is a random $N\times N$  matrix over \OFP\ with range contained in the set $\{M_\sigma\,|\,\sigma\in\S_N\}$. 

If $\S$ is a non-empty subset of $\S_N$,  a random $N\times N$ permutation matrix $U$ over \OFP\ is said to be \emph{uniformly distributed over $\S$} iff its range is the set $\{M_\sigma\,|\,\sigma\in\S\}$ and
\begin{equation}\label{uniform}
P\left[U^{-1}(M_\sigma)\right]=\frac{1}{\card\;\S}\qquad\forall\sigma\in\S.
\end{equation}
\end{random permutations}
\begin{notn}\label{notn1}
Throughout the remainder of this section, the following are fixed:

(a) A $*$-free family $(u_r)_{r=1}^s$ of random variables in a $*$-probability space $(\A,\varphi)$ such that each $u_r$ is a $d_r$-Haar unitary [notation as in \ref{notn s}(a)]. For every $w\in\F\setminus\{e\}$ [notation as in \ref{notn s}(b)], $u_w$ denotes the random variable in $(\A,\varphi)$ obtained from $w$ by replacing each $g_r$ by $u_r$ and each $g_r^*$ by $u_r^*$; by convention, $u_e:=1$ (the unit of the algebra \A).

It is immediate from the definitions that for every $w\in\F$, we have that
\begin{equation}\label{Haar word}\varphi(u_w)=
\begin{cases}1&\text{if $w\approx e$}
\\0&\text{if $w\not\approx e$.}
\end{cases}
\end{equation}

(b) A sequence $A_1,A_2,\ldots,A_s$ of (not necessarily disjoint, or even distinct) non-empty sets of positive integers such that $\sup A_r=d_r\;\forall r\in[s]$ and such that each $A_r$ satisfies condition \eqref{condition}, and an increasing sequence $(N_k)_{k=1}^\infty$ such that $\S_{N_k}^{(A_r)}\neq\emptyset\;\forall k\geq1,\forall r\in[s]$ (notation as in \ref{permutations}). For every graph $\Gamma$, $\mathcal C(\Gamma)$ denotes the set of all congruences $\pi\in\Con(\Gamma)$ such that $(\Gamma/\pi)(r)$ is an $A_r$-graph for every $r\in[s]$ [notation and terminology as in \ref{graphs}(a), \ref{quotient graphs}, \ref{admissible graphs}(b), and \ref{permutations and graphs}(a)].

(c) For every positive integer $k$, an independent family $\left(U_r\tok\right)_{r=1}^s$ of random $N_k\times N_k$ permutation matrices over \OFP\ such that each $U_r\tok$ is uniformly distributed over $\S_{N_k}^{(A_r)}$. The independence requirement implies that for every positive integer $k$ and for every measurable function $g:\left(\M_N\right)^s\to\C$, we have [with the notation from \ref{permutations} and \ref{random permutations}(a)] that
\begin{equation}\label{independent}
\int_\Omega g\left(U_1\tok(\omega),U_2\tok(\omega),\ldots,U_s\tok(\omega)\right)\dif P(\omega)=\frac{1}{\prod\limits_{r=1}^sa_{N_k}^{(A_r)}}\sum_{\substack{\sigma_1\in\S_{N_k}^{(A_1)}\\\sigma_2\in\S_{N_k}^{(A_2)}
\\\vdots
\\\sigma_s\in\S_{N_k}^{(A_s)}}}g(M_{\sigma_1},M_{\sigma_2},\ldots,M_{\sigma_s}).
\end{equation}
For every $w\in\F\setminus\{e\}$ [notation as in \ref{notn s}(b)] and for every positive integer $k$, $U_w\tok$ denotes the random $N_k\times N_k$ (permutation) matrix over \OFP\ obtained from $w$ by replacing each $g_r$ by $U_r\tok$ and each $g_r^*$ by $\left(U_r^{(k)}\right)^*$; by convention, $U_e\tok$ is the constant random matrix equal to $I\tok$ (the $N_k\times N_k$ identity matrix). 

If each of the sets $A_1,A_2,\ldots,A_s$ either consists of a single positive integer or is infinite, then it is easily seen that for every positive integer $k$ and for every $w\in\F$ such that $w\approx e$, we have that $U_w\tok$ is the constant random matrix equal to the $I\tok$.
\end{notn}
The main result of this paper is the following.
\begin{thm}\label{thm}
With the notation from \ref{notn1} and \ref{random matrices}, we have that:
\vspace{2mm}

(a) $\left(U_r^{(k)}\right)_{r=1}^s\converges(u_r)_{r=1}^s$.
\vspace{2mm}

(b) If each of the sets $A_1,A_2,\ldots,A_s$ either consists of a single positive integer or is infinite, then $\left(U_r^{(k)}\right)_{r=1}^s\convergesas(u_r)_{r=1}^s$.
\end{thm}
\begin{remks}
Theorem \ref{thm} is equivalent to the statement that
\begin{equation}\label{eqn1 thm1}
\lim_{k\to\infty}\E\left(\tr\;U_w\tok\right)=\varphi(u_w)\qquad\forall w\in\F\setminus\{e\}
\end{equation}
and that, if each of the sets $A_1,A_2,\ldots,A_s$ either consists of a single positive integer or is infinite, then
\begin{equation}\label{eqn2 thm1}
\lim_{k\to\infty}\tr\;U_w\tok=\varphi(u_w)\quad\text{almost surely}\qquad\forall w\in\F\setminus\{e\}.
\end{equation}
By Proposition \ref{prop as}, \eqref{eqn2 thm1} follows from \eqref{eqn1 thm1} provided that it is shown that, if each of the sets $A_1,A_2,\ldots,A_s$ either consists of a single positive integer or is infinite, then
\begin{equation}\label{eqn3 thm1}
\sum_{k=1}^\infty\left[\E\left(\left|\tr\;U_w\tok\right|^2\right)-\left|\E\left(\tr\;U_w\tok\right)\right|^2\right]<\infty
\qquad\forall w\in\F\setminus\{e\}.
\end{equation}
In turn, \eqref{eqn3 thm1} is a particular case of the more general statement that, if each of the sets $A_1,A_2,\ldots,A_s$ either consists of a single positive integer or is infinite, then
\begin{equation}\label{eqn4 thm1}
\sum_{k=1}^\infty\left[\E\left(\tr\;U_{w_1}\tok\cdot\tr\;U_{w_2}\tok\right)-\E\left(\tr\;U_{w_1}\tok\right)\cdot\E\left(\tr\;U_{w_2}\tok\right)\right]<\infty\qquad\forall w_1,w_2\in\F\setminus\{e\}.
\end{equation}
The following lemma gives an explicit formula for the type of expectations appearing in the above equations.
\end{remks}
\begin{lemma}\label{lemma thm}
With the notation from \ref{notn1}, \ref{graphs}(a), \ref{words}, \ref{quotient graphs}, and \ref{permutations and graphs}, let $n$ and $k$ be  positive integers and let $w_1,w_2,\ldots,w_n\in\F\setminus\{e\}$. Then
\begin{equation}\label{eqn2 lem2 thm1}
\E\left(\prod_{a=1}^n\tr\;U_{w_a}\tok\right)=\frac{1}{N_k^n}\sum_{\pi\in\mathcal C(\Gamma)}\frac{N_k!}{\big(N_k-|\pi|\big)!}\prod\limits_{r=1}^sp_{N_k}^{(A_r)}\big[(\Gamma/\pi)(r)\big],
\end{equation}
where $\Gamma$ denotes some (arbitrary) graph having $n$ distinct connected components $\Gamma_1,\Gamma_2,\ldots,\Gamma_n$ with $\Gamma_a\cong\Gamma_{w_a}\;\forall a\in[n]$. [The actual choice of $\Gamma$ is irrelevant, since the right-hand side of \eqref{eqn2 lem2 thm1} depends only on the isomorphism class of $\Gamma$.]
\begin{proof}
Write  $w_1=u_{r_1}^{\varepsilon_1}u_{r_2}^{\varepsilon_2}\cdots u_{r_{l_1}}^{\varepsilon_{l_1}},w_2=u_{r_{l_1+1}}^{\varepsilon_{l_1+1}}u_{r_{l_1+2}}^{\varepsilon_{l_1+2}}\cdots u_{r_{l_2}}^{\varepsilon_{l_2}},\ldots,w_n=u_{r_{l_{n-1}+1}}^{\varepsilon_{l_{n-1}+1}}u_{r_{l_{n-1}+2}}^{\varepsilon_{l_{n-1}+2}}\cdots u_{r_{l_n}}^{\varepsilon_{l_n}}$, where $l_1<l_2<\cdots<l_n$, $r_1,r_2,\ldots,r_{l_n}\in[s]$ and $\varepsilon_1,\varepsilon_2,\ldots,\varepsilon_{l_n}\in\{1,*\}$. Then for every positive integer $k$, we have (with the convention that $l_0:=0$) that
\begin{align*}
&\E\left(\prod_{a=1}^n\tr\;U_{w_a}\tok\right)=\E\left[\prod_{a=1}^n\tr\;\left(U_{r_{l_{a-1}+1}}\tok\right)^{\varepsilon
_{l_{a-1}+1}}\left(U_{r_{l_{a-1}+2}}\tok\right)^{\varepsilon_{l_{a-1}+2}}\cdots\left(U_{r_{l_a}}\tok\right)^{\varepsilon_{l_a}}\right]
\\&=\int_\Omega\frac{1}{N_k^n}\sum_{i:[l_n]\to[N_k]}\prod_{a=1}^n\bigg[\left(U_{r_{l_{a-1}+1}}\tok(\omega)^{\varepsilon_{l_{a-1}+1}}\right)_{i(l_{a-1}+1),i(l_{a-1}+2)}\cdot
\\&\hspace{3cm}\cdot\left(U_{r_{l_{a-1}+2}}\tok(\omega)^{\varepsilon_{l_{a-1}+2}}\right)_{i(l_{a-1}+2),i(l_{a-1}+3)}\cdots\left(U_{r_{l_a}}\tok(\omega)^{\varepsilon_{l_a}}\right)_{i(l_a)
,i(l_{a-1}+1)}\bigg]\dif P(\omega)
\\&=\frac{1}{N_k^n}\sum_{i:[l_n]\to[N_k]}\frac{1}{\prod\limits_{r=1}^sa_{N_k}^{(A_r)}}\sum_{\substack{\sigma_1\in\S_{N_k}^{(A_1)}
\\\sigma_2\in\S_{N_k}^{(A_2)}
\\\vdots
\\\sigma_s\in\S_{N_k}^{(A_s)}}}\prod_{a=1}^n\bigg[\left(M_{\sigma_{r_{l_{a-1}+1}}}^{\varepsilon_{l_{a-1}+1}}\right)_{i({l_{a-1}+1}),i({l_{a-1}+2})}\cdot
\\&\hspace{3cm}\cdot\left(M_{\sigma_{r_{l_{a-1}+2}}}^{\varepsilon_{l_{a-1}+2}}\right)_{i({l_{a-1}+2}),i({l_{a-1}+3})}\cdots\left(M_{\sigma_{r_{l_a}}}^{\varepsilon_{l_a}}\right)_{i(l_a),i(l_{a-1}+1)}\bigg]\tag*{\text{[by \eqref{independent}]}}
\\&=\frac{1}{N_k^n}\sum_{i:[l_n]\to[N_k]}\frac{1}{\prod\limits_{r=1}^sa_{N_k}^{(A_r)}}\;\card\left\{(\sigma_1,\sigma_2,\ldots,\sigma_s)\,\left|\,\substack{\sigma_r\in\S_{N_k}^{(A_r)}\;\forall r\in[s]\text{ and}
\\\sigma_{r_a}^{\varepsilon_a}\big(i[\tau(l)]\big)=i(l)\;\forall l\in[l_n]}\right\}\right.,
\end{align*}
with the convention that $\sigma^*:=\sigma^{-1}\;\forall\sigma\in\S_{N_k}$, and where $\tau\in\S_{l_n}$ is the permutation defined by
$$\tau=(1,2,\ldots,l_1)(l_1+1,l_1+2,\ldots,l_2)\cdots(l_{n-1}+1,l_{n-1}+2,\ldots,l_n).$$
Next note that every function $i:[l_n]\to [N_k]$ is determined by a pair $(\pi,f)$, where $\pi\in\Pi\big([l_n]\big)$ is defined by
$$a\overset{\pi}{\sim}b\iff i(a)=i(b)\qquad\forall a,b\in [l_n]$$
and $f:\pi\to [N_k]$ is the (injective) function defined by
$$f\big(\pi(a)\big)=i(a)\qquad\forall a\in [l_n].$$
Let $\Gamma$ be the graph having the $n$ distinct connected components $\Gamma_{w_1},\Gamma_{w_2},\ldots,\Gamma_{w_n}$, with each $\Gamma_{w_a}$ having the vertices $1,2,\ldots,l_a-l_{a-1}$ relabeled as $l_{a-1}+1,l_{a-1}+2,\ldots,l_a$. Then
\begin{align*}\E\left(\prod_{a=1}^n\tr\;U_{w_a}\tok\right)&=\frac{1}{N_k^n}\sum_{\substack{\pi\in\Pi\big([l_n]\big)
\\f:\pi\to[N_k]
\\\text{injective}}}\frac{1}{\prod\limits_{r=1}^sa_{N_k}^{(A_r)}}\;\card\left\{(\sigma_1,\sigma_2,\ldots,\sigma_s)\,\left|\,\substack{\sigma_r\in\S_{N_k}^{(A_r)}\;\forall r\in[s]\text{ and}
\\\sigma_{r_a}^{\varepsilon_a}\big[f\big(\pi[\tau(l)]\big)\big]=f\big(\pi(l)\big)\;\forall l\in[l_n]}\right\}\right.
\\&=\frac{1}{N_k^n}\sum_{\pi\in\Pi\big([l_n]\big)}\sum_{\substack{f:\pi\to[N_k]
\\\text{injective}}}\frac{1}{\prod\limits_{r=1}^sa_{N_k}^{(A_r)}}\;\prod\limits_{r=1}^s\card\;\S_{N_k}^{(A_r)}\big[f(\Gamma/\pi)(r)\big]\tag*{\text{[notation as in \ref{permutations and graphs}]}}
\\&=\frac{1}{N_k^n}\sum_{\pi\in\mathcal C(\Gamma)}\prod_{r=1}^sp_{N_k}^{(A_r)}\big[(\Gamma/\pi)(r)\big]\cdot\card\;\big\{f:\pi\to[N_k]\,\big|\,f\text{ injective}\big\}
\\&=\frac{1}{N_k^n}\sum_{\pi\in\mathcal C(\Gamma)}\frac{N_k!}{\big(N_k-|\pi|\big)!}\prod\limits_{r=1}^sp_{N_k}^{(A_r)}\big[(\Gamma/\pi)(r)\big].
\end{align*}
\end{proof}
\end{lemma}
\pagebreak

\begin{coro}\label{coro thm}
With the notation from \ref{notn1}, we have that:

(a) For every $w_1,w_2,\ldots,w_n\in\F$, 
\begin{equation}\label{coro1 thm1}
\lim_{k\to\infty}\E\left(\prod_{a=1}^n\tr\;U_{w_a}\tok\right)=
\begin{cases}1&\text{if $w_a\approx e\;\forall a\in[n]$}
\\0&\text{otherwise.}
\end{cases}
\end{equation}
In addition, if each of the sets $A_1,A_2,\ldots,A_s$ either consists of a single positive integer or is infinite and it is not the case that $w_a\approx e\;\forall a\in[n]$, then
\begin{equation}\label{coro2 thm1}
\E\left(\prod_{a=1}^n\tr\;U_{w_a}\tok\right)=\on\qquad\text{as $k\to\infty$}.
\end{equation}

(b) If each of the sets $A_1,A_2,\ldots,A_s$ either consists of a single positive integer or is infinite, then for every $w_1,w_2\in\F$ with $w_1\not\approx e\not\approx w_2$, there exists $\delta>0$ such that
\begin{equation}\label{coro3 thm1}
\E\left(\tr\;U_{w_1}\tok\cdot\tr\;U_{w_2}\tok\right)=\onn\qquad\text{as $k\to\infty$}.
\end{equation}
\begin{proof}
(a) It may be supposed that $w_a\neq e\;\forall a\in[n]$, since for every positive integer $k$, we have that $U_e\tok$ is the constant random matrix equal to the $N_k\times N_k$ identity. Let $\Gamma$ be an arbitrary graph with $n$ connected components $\Gamma_1,\Gamma_2,\ldots,\Gamma_n$, where each $\Gamma_a\cong\Gamma_{w_a}$. With the notation from \ref{notn1},
\ref{graphs}(a), \ref{words}, \ref{quotient graphs}, \ref{loop-characteristic}, and \ref{permutations and graphs}, we have that
\begin{align*}
\E\left(\prod_{a=1}^n\tr\;U_{w_a}\tok\right)&\sim\sum_{\pi\in\mathcal C(\Gamma)}N_k^{|\pi|-n}\prod\limits_{r=1}^sp_{N_k}^{(A_r)}\big[(\Gamma/\pi)(r)\big]\tag*{\text{[by Lemma \ref{lemma thm}]}}
\\&\sim\sum_{\pi\in\mathcal C(\Gamma)}N_k^{V(\Gamma/\pi)-n}\prod\limits_{r=1}^sN_k^{-E\big[(\Gamma/\pi)(r)\big]+L\big[(\Gamma/\pi)(r)\big]+\sum\limits_{\substack{\text{loops $P$ of}
\\\text{$(\Gamma/\pi)(r)$}}}\left[\frac{l(P)}{d_r}-1\right]}\tag*{\text{[by Proposition \ref{limpng}]}}
\\&\sim\sum_{\pi\in\mathcal C(\Gamma)}N_k^{\chi(\Gamma/\pi)-n+\sum\limits_{r=1}^s\sum\limits_{\substack{\text{$r$-loops $P$}
\\\text{of $\Gamma/\pi$}}}\left[\frac{l(P)}{d_r}-1\right]}\qquad\text{as $k\to\infty$}.
\end{align*}
For every $\pi\in\mathcal C(\Gamma)$, we have that 
$\chi(\Gamma/\pi)-n+\sum\limits_{r=1}^s\sum\limits_{\substack{\text{$r$-loops $P$}
\\\text{of $\Gamma/\pi$}}}\left[\frac{l(P)}{d_r}-1\right]\leq0$
by Proposition \ref{connected chi} and the definition of $\mathcal C(\Gamma)$, with equality holding iff $\pi\in\SCon(\Gamma)$ (notation as in \ref{strongly admissible graphs}) and $\Gamma/\pi$ has $n$ connected components, each having loop-characteristic 1; that is, iff $\pi$ is of the form $\bigcup\limits_{a=1}^n\pi_a$, where each $\pi_a\in\SCon(\Gamma_a)$ and $\chi(\Gamma_a/\pi_a)=1$. Thus
\begin{equation*}
\lim_{k\to\infty}\E\left(\prod_{a=1}^n\tr\;U_{w_a}\tok\right)=\prod_{a=1}^n\card\;\big\{\pi\in\SCon(\Gamma_a)\,\big|\,\chi(\Gamma_a/\pi)=1\big\}=
\begin{cases}1&\text{if $w_a\approx e\;\forall a\in[n]$}
\\0&\text{otherwise.}
\end{cases}\tag*{\text{[by Proposition \ref{prop scon chi word-graphs}]}}
\end{equation*}
In addition, if each of the sets $A_1,A_2,\ldots,A_s$ either consists of a single positive integer or is infinite and it is not the case that $w_a\approx e\;\forall a\in[n]$, then $\chi(\Gamma/\pi)-n+\sum\limits_{r=1}^s\sum\limits_{\substack{\text{$r$-loops $P$}\\\text{of $\Gamma/\pi$}}}\left[\frac{l(P)}{d_r}-1\right]$ is a negative integer for every $\pi\in\mathcal C(\Gamma)$, so it follows that
$$\E\left(\prod_{a=1}^n\tr\;U_{w_a}\tok\right)=\on\qquad\text{as $k\to\infty$}.$$

(b) Suppose that each of the sets $A_1,A_2,\ldots,A_s$ either consists of a single positive integer or is infinite, and let $w_1,w_2\in\F$ be such that $w_1\not\approx e\not\approx w_2$. Let $\Gamma$ be an arbitrary graph with two connected components $\Gamma_1$ and $\Gamma_2$, where $\Gamma_1\cong\Gamma_{w_1}$ and $\Gamma_2\cong\Gamma_{w_2}$. Exactly as in the
proof of part (a), we have that
$$\E\left(\tr\;U_{w_1}\tok\cdot\tr\;U_{w_2}\tok\right)\sim\sum_{\pi\in\mathcal C(\Gamma)}N_k^{\chi(\Gamma/\pi)-2+\sum\limits_{r=1}^s\sum\limits_{\substack{\text{$r$-loops $P$}
\\\text{of $\Gamma/\pi$}}}\left[\frac{l(P)}{d_r}-1\right]}\qquad\text{as $k\to\infty$},$$
so it suffices to show that $\chi(\Gamma/\pi)-2+\sum\limits_{r=1}^s\sum\limits_{\substack{\text{$r$-loops $P$}\\\text{of $\Gamma/\pi$}}}\left[\frac{l(P)}{d_r}-1\right]<-1$ for every $\pi\in\mathcal C(\Gamma).$

Suppose first that $\pi\in\mathcal C(\Gamma)$ is of the form $\pi_1\cup\pi_2$ for some partitions $\pi_1$ and $\pi_2$ of the vertex sets of $\Gamma_1$ and $\Gamma_2$, respectively. Then $\Gamma/\pi$ consists of the two connected components $\Gamma_1/\pi_1$ and $\Gamma_2/\pi_2$, so
\begin{align*}
&\chi(\Gamma/\pi)-2+\sum\limits_{r=1}^s\sum\limits_{\substack{\text{$r$-loops $P$}
\\\text{of $\Gamma/\pi$}}}\left[\frac{l(P)}{d_r}-1\right]\\&=\left[\chi(\Gamma_1/\pi_1)-1+\sum\limits_{r=1}^s\sum\limits_{\substack{\text{$r$-loops $P$}
\\\text{of $\Gamma_1/\pi_1$}}}\left[\frac{l(P)}{d_r}-1\right]\right]+\left[\chi(\Gamma_2/\pi_2)-1+\sum\limits_{r=1}^s\sum\limits_{\substack{\text{$r$-loops $P$}
\\\text{of $\Gamma_2/\pi_2$}}}\left[\frac{l(P)}{d_r}-1\right]\right].
\end{align*}
As in the proof of part (a), each of the two summands is a negative integer, so it follows that $\chi(\Gamma/\pi)-2+\sum\limits_{r=1}^s\sum\limits_{\substack{\text{$r$-loops $P$}\\\text{of $\Gamma/\pi$}}}\left[\frac{l(P)}{d_r}-1\right]\leq-2$.

Finally, suppose that $\pi\in\mathcal C(\Gamma)$ is not of the form $\pi_1\cup\pi_2$ for some partitions $\pi_1$ and $\pi_2$ of the vertex sets of $\Gamma_1$ and $\Gamma_2$, respectively. Then $\Gamma/\pi$ is connected, so to show that $\chi(\Gamma/\pi)-2+\sum\limits_{r=1}^s\sum\limits_{\substack{\text{$r$-loops $P$}\\\text{of $\Gamma/\pi$}}}\left[\frac{l(P)}{d_r}-1\right]<-1$, it suffices to show that either $\Gamma/\pi$ is not strongly admissible or $\chi(\Gamma/\pi)<1$. Suppose that $\Gamma/\pi$ is strongly admissible and that $\chi(\Gamma/\pi)=1$, and let $\pi_1$ be the restriction of $\pi$ to the vertices of $\Gamma_1$, as in \ref{partitions}(a). Then $\Gamma_1/\pi_1$ is isomorphic to a subgraph of $\Gamma/\pi$, so $\Gamma_1/\pi_1$ is strongly admissible and $\chi(\Gamma_1/\pi_1)=1$ by Remarks \ref{subgraphs}. Thus $\pi_1\in\SCon(\Gamma_1)$, which contradicts Proposition \ref{prop scon chi word-graphs} (since $w_1\not\approx e$).
\end{proof}
\end{coro}
\textit{Proof of Theorem \ref{thm}.}
(a) It suffices to prove \eqref{eqn1 thm1}, and this follows directly from \eqref{Haar word} and Corollary \ref{coro thm}(a).

(b) It suffices to prove \eqref{eqn4 thm1}. Suppose that each of the sets $A_1,A_2,\ldots,A_s$ either consists of a single positive integer or is infinite, and let $w_1,w_2\in\F\setminus\{e\}$.

If $w_1\approx e$, then for every positive integer $k$, we have that $U_{w_1}\tok$ is the constant random matrix equal to 
the $N_k\times N_k$ identity [see the remark at the end of \ref{notn1}(c)], so
$$\E\left(\tr\;U_{w_1}\tok\cdot\tr\;U_{w_2}\tok\right)-\E\left(\tr\;U_{w_1}\tok\right)\cdot\E\left(\tr\;U_{w_2}\tok\right)=\E\left(\tr\;U_{w_2}\tok\right)-\E\left(\tr\;U_{w_2}\tok\right)=0.$$
Similarly, if $w_2\approx e$, then
$\E\left(\tr\;U_{w_1}\tok\cdot\tr\;U_{w_2}\tok\right)-\E\left(\tr\;U_{w_1}\tok\right)\cdot\E\left(\tr\;U_{w_2}\tok\right)=0.$

Finally, if $w_1\not\approx e\not\approx w_2$, then
$$\E\left(\tr\;U_{w_1}\tok\right)\cdot\E\left(\tr\;U_{w_2}\tok\right)=\mathrm O\left(1/N_k^2\right)\qquad\text{as $k\to\infty$}$$
and there exists $\delta>0$ such that
$$\E\left(\tr\;U_{w_1}\tok\cdot\tr\;U_{w_2}\tok\right)=\onn\qquad\text{as $k\to\infty$}$$
by Corollary \ref{coro thm}, so
$\sum\limits_{k=1}^\infty\left[\E\left(\tr\;U_{w_1}\tok\cdot\tr\;U_{w_2}\tok\right)-\E\left(\tr\;U_{w_1}\tok\right)\cdot\E\left(\tr\;U_{w_2}\tok\right)\right]<\infty.$\qed
\begin{gaussian matrices}\label{gaussian matrices}
Within the framework of Theorem \ref{thm}, the family $\left(U_r\tok\right)_{r=1}^s$ is also asymptotically $*$-free from either an independent family of square complex Gaussian matrices, or from an independent family of complex Wishart matrices. More concretely, let $(f_{r;i,j})_{\substack{r\in[s]\\i,j\geq1}}$ be an independent family of complex standard Gaussian random variables over $(\Omega,\mathcal F,P)$; in other words, $(\mathrm{Re}\;f_{r;i,j},\mathrm{Im}\;f_{r;i,j})_{\substack{r\in[s]\\i,j\geq1}}$ is an independent family of real Gaussian random variables over $(\Omega,\mathcal F,P)$ with mean 0 and variance 1/2. For every $M,N\geq1$, define random $M\times N$ matrices $G_1\toMN,G_2\toMN,\ldots,G_s\toMN$ by
\begin{equation}
\left(G_r\toMN\right)_{i,j}=f_{r;i,j}\qquad\forall i\in[M],\forall j\in[N].
\end{equation}
Fix a sequence $(M_k)_{k=1}^\infty$ of integers with the ratio $M_k/N_k$ converging to some real number $c$ and for every $r\in[s]$ and $k\geq1$, let $G_r^{(k)}:=\frac{1}{\sqrt N_k}G_r^{(N_k,N_k)}$ and $W_r^{(k)}:=\frac{1}{N}\left(G_r^{(M_k,N_k)}\right)^*G_r^{(M_k,N_k)}$. (The random $N_k\times N_k$ matrices $W_r^{(k)}$ are known as \emph{complex Wishart
matrices}.) Then an almost identical proof to that of Theorem 3.1 of \cite{Neagu} shows that:

(1) If the family $\left(U_r^{(k)},G_r^{(k)}\right)_{r=1}^s$ is independent for every $k\geq1$, then
$\left(U_r^{(k)},G_r^{(k)}\right)_{r=1}^s\converges(u_r,g_r)_{r=1}^s,$
where each $g_r$ is a \emph{standard circular} random variable in $(\mathcal A,\varphi)$ (see \cite{Neagu} for the definition) and the family $(u_r,g_r)_{r=1}^s$ is $*$-free.

(2) If the family $\left(U_r^{(k)},W_r^{(k)}\right)_{r=1}^s$ is independent for every $k\geq1$, then
$\left(U_r^{(k)},W_r^{(k)}\right)_{r=1}^s\converges(u_r,w_r)_{r=1}^s,$
where each $w_r$ is a \emph{free Poisson} random variable of parameter $c$ in $(\mathcal A,\varphi)$ (see \cite{Neagu} for the definition) and the family $(u_r,w_r)_{r=1}^s$ is $*$-free.

In both cases, if each of the sets $A_1,A_2,\ldots,A_s$ either consists of a single positive integer or is infinite, then the convergence in $*$-distribution actually holds almost surely.
\end{gaussian matrices}
\bibliographystyle{amsplain}
\bibliography{refs}

\providecommand{\bysame}{\leavevmode\hbox to3em{\hrulefill}\thinspace}
\begin{thebibliography}{1}

\bibitem{Bona}
M.~B\'ona, \emph{Combinatorics of {Permutations}}, Chapman \& Hall/CRC, Boca
  Raton, FL, 2004.

\bibitem{Hildebrand}
A.~Hildebrand and G.~Tenenbaum, \emph{On some {Tauberian} theorems related to
  the prime number theorem}, Compositio Math. \textbf{90} (1994), no.~3,
  315--349.

\bibitem{Manstavicius}
E.~{Manstavi\v cius}, \emph{Mappings on decomposable combinatorial structures:
  {Analytic} approach}, Combin. Probab. Comput. \textbf{11} (2002), no.~1,
  61--78.

\bibitem{Neagu}
M.~G. Neagu, \emph{Asymptotic freeness of random permutation matrices from
  {Gaussian} matrices}, J. Ramanujan Math. Soc. \textbf{20} (2005), no.~3,
  189--213.

\bibitem{Nica1}
A.~Nica, \emph{Asymptotically free families of random unitaries in symmetric
  groups}, Pacific J. Math. \textbf{157} (1993), 295--310.

\bibitem{Nica2}
\bysame, \emph{On the number of cycles of given length of a free word in
  several random permutations}, Random Structures Algorithms \textbf{5} (1994),
  703--730.

\bibitem{Thorbjornsen}
S.~Thorbj\o rnsen, \emph{Mixed moments of {Voiculescu's} {Gaussian} random
  matrices}, J. Funct. Anal. \textbf{176} (2000), 213--246.

\bibitem{Voiculescu}
D.~Voiculescu, \emph{Limit laws for random matrices and free products}, Invent.
  Math. \textbf{104} (1991), 201--220.

\bibitem{Wilf}
H.~Wilf, \emph{{Generatingfunctionology}}, 2nd ed., Academic Press, Inc.,
  Boston, MA, 1994.

\end{thebibliography}
\end{document}